\documentclass[9pt,shortpaper,twoside,web]{ieeecolor}

\pdfoutput=1 
\usepackage{hyperref}
\usepackage{generic}
\usepackage{cite}
\usepackage{amsmath,amssymb,amsfonts}
\usepackage{algorithmic}
\usepackage{graphicx}

\usepackage{mathtools} 
\usepackage[most]{tcolorbox} 
\usepackage{caption} 
\usepackage{subcaption}

\newtheorem{thm}{Theorem}
\newtheorem{lemma}{Lemma}

\newtheorem{remark}{Remark}
\newtheorem{asmp}{Assumption} 

\newcommand{\reals}{\mathbb{R}}

\newcommand{\norm}[1]{\left \lVert#1\right \rVert}
\newcommand{\set}[1]{\{#1\}}
\newcommand{\setc}[2]{\{#1\ |\ #2\}}

\newcommand{\orth}{\perp}
\newcommand{\orthc}{||}

\newcommand{\Id}{\mathrm{Id}}

\DeclareMathOperator*{\argmin}{arg\,min}

\DeclareMathOperator{\diag}{diag}

\DeclarePairedDelimiterX{\inp}[2]{\langle}{\rangle}{#1\ \middle| \ #2}

\usepackage{multirow}
\begin{document}
\title{Single-timescale distributed GNE seeking for  aggregative games  over networks via forward-backward operator splitting} 

\author{Dian Gadjov and Lacra Pavel
\thanks{This work was supported by NSERC Discovery Grant.}%
\thanks{\footnotesize D. Gadjov and L. Pavel are with Dept. of Electrical and Computer Engineering, University of Toronto. 
        {\tt\footnotesize dian.gadjov@mail.utoronto.ca}, 
        {\tt\footnotesize pavel@control.utoronto.ca}}%
}

\maketitle

\begin{abstract}
We consider  aggregative games with affine coupling constraints, where agents have partial information on the aggregate value and can only communicate with neighbouring agents. We propose a single-layer distributed algorithm that reaches a variational generalized Nash equilibrium, under constant step sizes. 
The algorithm works on a single timescale, i.e., does not require multiple communication rounds between agents before updating their action. The convergence proof leverages an invariance property of the aggregate estimates and relies on a forward-backward splitting for two preconditioned operators  and  their restricted (strong) monotonicity properties on the consensus subspace. 
\vspace{-0.3cm}\end{abstract}
\vspace{-0.3cm}
\section{Introduction}

Research in aggregative games has surged in recent years, due to their suitability to model decision problems in various application domains: from  demand-side management for the smart grid, \cite{Basar2012} and electric vehicles, 
 \cite{PEVMAAuto},   to 
network congestion control, \cite{Garcia}. 
Aggregative games are non-cooperative games in which each (player's) agent's cost depends on some aggregate effect of all other agents' actions. Often agents have shared coupling constraints, 
and the relevant equilibrium concept is the generalized Nash equilibrium (GNE).

Many settings involve a large number of agents with private cost functions and constraints, who are willing  and able to exchange information with  
their neighbours only, \cite{Tekin2012}, \cite{cournotgame}. 
Motivated by the above,  in this note we develop a single-layer/single-timescale, { discrete-time} distributed GNE seeking algorithm for \emph{aggregative games}, which is guaranteed to converge exactly to a variational GNE when using fixed-step sizes. Agents update simultaneously their action and an aggregate-estimate, based on local communication. To the best of our knowledge, this is the first such algorithm in the literature {for aggregative games}. Our novel contributions are to relate the algorithm to a preconditioned forward-backward splitting of a pair of monotone operators,  and develop conditions for distributed convergence.

\emph{Literature Review:} The problem of finding a (generalized) Nash Equilibrium (G)NE when agents know the actions of all other agents (full-decision information) has been studied thoroughly, e.g.  \cite{arrow,B99,FacchineiGNEP,FP07}, or via a centralized operator-splitting \cite{combettesBriceno}. In recent years, a rapidly growing field  originated by \cite{NedichDistAgg,FPAuto2016} is focused on developing algorithms that relax this full-information assumption to (G)NE computation over networks. {Fully distributed  NE  seeking algorithms have been proposed for general games with decoupled constraints, either  in continuous time, e.g. \cite{cat},\cite{Gadjov} 
or in discrete time (with fixed-step sizes), e.g. \cite{SalehiShiPavelAuto_2019}, or \cite{TatarenkoShi} (with fast convergence).  However, these algorithms require each agent to maintain an estimate of the actions of all other agents,  hence are  inefficient for aggregative games. The same holds for the fully distributed GNE seeking algorithm proposed  in  \cite{Lacra} for generally coupled games.} 
For \emph{aggregative games}, existing algorithms can be classified  as \emph{semi-decentralized} (when a central coordinator is  required), or  \emph{distributed} (when only local communication is used). \emph{Semi-decentralized GNE algorithms}  have been proposed  for \emph{aggregative games} in \cite{GrammaticoTAC2017,belg,wardrop,gramDR} (discrete-time) and \cite{PersisGrammatico}  (continuous-time). An elegant operator theory approach is used to show that global convergence can be achieved with fixed-step sizes under (strict) strong monotonicity of the pseudo-gradient of the game, either to a variational generalized Nash equilibrium (GNE),  or to an aggregative (Wardrop) equilibrium (GAE). 
The algorithms require  a central coordinator to broadcast the aggregate value and to ensure the coupling constraints are met. 
  
In a distributed setting a central coordinator/node does not exist, and distributed algorithms are more difficult to develop. For aggregative games, such distributed algorithms have been mostly developed for NE seeking with no coupled constraints. The first work to propose \emph{distributed NE algorithms} for \emph{aggregative games over networks} was \cite{NedichDistAgg}. In order to cope with not knowing the aggregate value, each agent maintains an estimate of the true aggregate,  built based on local communication with neighbours, and uses it in the action update  instead of  the true aggregate. 
The algorithm proposed in \cite{NedichDistAgg} requires diminishing step-sizes for exact convergence, while only guaranteeing convergence to a neighbourhood of the NE when using fixed-step sizes, under (strict) strong monotonicity of the pseudo-gradient of the game.  The recent algorithm proposed in \cite{Shanbhag} requires an increasing number of communication rounds before each action update. In both cases, the NE seeking algorithms 
effectively operate as if on two time-scales (fast aggregate estimate, slow action update) with each agent using an estimate that is near the true aggregate. 

Fewer results exist for distributed GNE seeking in aggregative games with \emph{coupled constraints}. In fact, the only such distributed discrete-time algorithm {{for aggregative games}}  that we are aware of is the one proposed  in \cite{parise}. The algorithm requires that players exchange information for a fixed number of communication rounds before every action update and is only guaranteed to reach an $\epsilon$-GNE. The $\epsilon$-GNE approaches a variational GNE if the fixed number of communication rounds between action updates goes to infinity, which is impractical. A \emph{continuous-time} dynamics for  aggregative games with equality constraints is proposed in \cite{PengAgg}.  However, 
convergence in continuous-time does not guarantee convergence in  discrete-time. 

\emph{Contributions: } Motivated by the above, in this paper we propose  a  distributed, single-layer discrete-time GNE seeking algorithm for aggregative games with affine coupling constraints, that has guaranteed exact convergence { to} a variational GNE when using constant step-sizes. To the best of our knowledge, to date there does not exist such an algorithm  { tailored for \emph{aggregative games}}. { We note that the  algorithm recently proposed in \cite{Lacra} for generally coupled games is applicable to aggregative games, but because it does not exploit the aggregative structure is computationally inefficient when used for these types of games.}
In our algorithm, each agent maintains an estimate of the aggregate and multiplier, and exchanges them with its neighbours over an undirected static connected graph, in order to learn the true aggregate value and ensure that the coupling constraints are satisfied.  Each agent updates its action  and its aggregate estimates in the same iteration. 
Compared to the distributed NE algorithm in \cite{NedichDistAgg}, and the GNE algorithm in  \cite{parise}  for \emph{aggregative games}, our algorithm does not require diminishing step-sizes nor multiple communication rounds between each action update. 
It  
is related to a preconditioned forward-backward operator splitting iteration, inspired by the distributed  framework for GNE seeking in general games conceptualized in \cite{PengLacraAuto2019}  for full-decision information, and extended  in \cite{Lacra} to partial-decision information. We note that 
the algorithm in \cite{Lacra} requires each agent to keep an estimate of \emph{all}  other agents' actions;  if  applied to \emph{aggregative games} where the coupling is only through the average, this can be not scalable/inefficient. 

The distributed algorithm we develop here 
exploits the aggregative coupling structure in the cost.  Each agent exchanges and maintains only an estimate of  the aggregate (dimensionally independent on the number of agents, hence scalable) and an estimate of the dual multiplier. Thus players do not need to share action information, which might be private information. 
We use proof techniques  similar  to  those in  \cite{Lacra} with the following differences. 
{ The algorithm is tailored for the special structure of aggregative games. 
Unlike  \cite{Lacra}, 
 here the estimate is a separate variable that needs to track the aggregate decision. Because of this, the update equation of  the aggregate estimate has an additional correction term. This is introduced to account for the own action's effect on the average and acts as a dynamic-tracking term. \emph{Secondly}, this correction term allows us to exploit an invariance property of the aggregate estimate. This invariance property plays a critical role in ensuring convergence  with fixed-step sizes and in obtaining a better bound compared to the algorithm in  \cite{Lacra}. \emph{Thirdly}, because of this correction term we need to use slightly different splitting operators and a different metric matrix, which needs handled separately. }  A conference version appears in \cite{DianLP2019}; here we provide   all proofs and additional numerical results. 

\emph{Notations}. For a vector $x\in\reals^{m}$, $x^{T}$ denotes its transpose and $\norm{x} = \sqrt{\inp*{x}{x}} = \sqrt{x^{T}x}$ the norm induced by inner product $\inp*{\cdot}{\cdot}$. For a symmetric positive-definite matrix $\Phi$, $\Phi \succ 0$, $\lambda_{min}(\Phi)$ and $\lambda_{max}(\Phi)$ denote its minimum and maximum eigenvalues. The $\Phi$-induced inner product is $\inp*{x}{y}_{\Phi} = \inp*{\Phi x}{y}$ and the $\Phi$-induced norm, $\norm{x}_{\Phi} = \sqrt{\inp*{\Phi x}{x}}$. For a matrix $A\in\reals^{m\times n}$, let $\norm{A} = \sigma_{max}(A)$ denote the 2-induced matrix norm, where $\sigma_{max}(A)$ is its maximum singular value. For $\mathcal{N} = \set{1, \dots, N}$, $col(x_{i})_{i\in \mathcal{N}}$ denotes the stacked vector obtained from vectors $x_{i}$, $diag(x_{i})_{i\in\mathcal{N}}$ is the diagonal matrix with $x_{i}$ along the diagonal. $Null(A)$ and $Range(A)$ are the null and range space of matrix $A$, while $[A]_{ij}$ stands for its $(i, j)$ entry. {  $I_{n}$, $\mathbf{1}_{n}$ and $\mathbf{0}_{n}$ denote the identity matrix, the all-ones and the all-zero vector of dimension $n$, respectively. We may also simply use $\mathbf{0}$ to denote an all-zero matrix of appropriate dimensions. } Denote $\prod_{i=1}^{N} \Omega_{i}$ as the Cartesian product of  sets $\Omega_{i}$, $i=1,\dots,N$. For a function $f(x) = f(col(x_{i})_{i\in \mathcal{N}})$ 
let $\nabla_{x_{i}} f(x) =  \frac{\partial}{\partial x_{i}}f(x)$. 

\vspace{-0.3cm}
\section{Background}\label{sec:background}
\vspace{-0.3cm}
\subsection*{Monotone Operators}
The following are from \cite{monoBookv1}. Let $\mathcal{A}: \reals^{m} \to 2^{\reals^{m}}$ be a set-valued operator. The domain of $\mathcal{A}$ is $dom \mathcal{A} = \setc{x \in \reals^{m}}{\mathcal{A}x \neq \emptyset}$ where $\emptyset$ is the empty set, and the range of $\mathcal{A}$ is $ran \mathcal{A} = \setc{y\in\reals^{m}}{ \exists x, y \in \mathcal{A}x}$. The graph of $\mathcal{A}$ is $gra \mathcal{A} = \setc{(x,u)\in\reals^{m}\times \reals^{m}}{u\in \mathcal{A}x}$. The zero set of $\mathcal{A}$ is $zer \mathcal{A} = \setc{x\in\reals^{m}}{\mathbf{0} \in \mathcal{A}x}$.  $\mathcal{A}$ is called monotone if $\forall (x, u), \forall (y, v) \in gra \mathcal{A}$, $\inp*{x-y}{u-v} \geq 0$, and strongly monotone if $\exists \mu > 0$ such that $\forall (x, u), \forall (y, v) \in gra \mathcal{A}$, $\inp*{x-y}{u-v} \geq \mu \norm{x-y}^{2}$. $ \mathcal{A}$ is maximally monotone if $gra \mathcal{A}$ is not strictly contained in the graph of any other monotone operator. The resolvent of $\mathcal{A}$ is $\mathcal{J}_{\mathcal{A}} = (\Id + \mathcal{A})^{-1}$, where $\Id$ is the identity operator. The fixed points of $\mathcal{A}$ are $Fix \mathcal{A} = \setc{x\in\reals^{m}}{x\in \mathcal{A}x}$.

An operator $T \!:\! \Omega \!\subset \! \reals^{m}\!\to \! \reals^{m}$ is nonexpansive if $\norm{T(x)\! -\! T(y)}  \!\leq  \!\norm{x \!-\! y},$ $\! \forall x, \!y \!\in \!\Omega$. $T \! \in \! \mathfrak{A}(\alpha)$, where $\mathfrak{A}(\alpha)$ denotes the class of $\alpha$-averaged operators, if and only if $\forall x, y \! \in\!\Omega$, $\norm{Tx\!-\!Ty}^{2} \!\leq \! \norm{x\!-\!y}^{2} \!-\! \frac{1\!-\!\alpha}{\alpha}\norm{(x-y)\!-\!(Tx-Ty)}^{2}$. $T\!\in\!\mathfrak{A}(\frac{1}{2})$ is also called firmly nonexpansive. If $\mathcal{A}$ is maximally monotone then $\mathcal{J}_{\mathcal{A}}$ is firmly nonexpasive,
[\cite{monoBookv1}, Prop. 23.7]. Let the projection of $x$ onto $\Omega$ be $P_{\Omega}(x) \!=\! \argmin_{y\in\Omega} \norm{x \!-\! y}$, with $P_{\Omega}(x)
\!=\! \mathcal{J}_{N_{\Omega}}(x)\!=(\Id \!+\! N_{\Omega})^{-1}$, where  $N_{\Omega}(x) \!=\! \setc{v}{\inp*{v}{y \!-\! x} \!\leq \!0,\! \forall y\!\in\!\Omega}$ is the normal cone.   
If $\Omega$ is closed and convex, $P_{\Omega}$ is firmly nonexpansive 
[\cite{monoBookv1}, Prop. 4.8]. $T$ is called $\beta$-cocoercive if $\beta T \! \in \!\mathfrak{A}(\frac{1}{2})$ for $\beta \!>\! 0$, i.e., $\!\beta\norm{Tx \!- Ty}^{2}\! \leq \! \inp*{x-y}{Tx\!-\!Ty}$ $\forall x, y \!\in \!\Omega$.\!
\vspace{-0.25cm}
\subsection*{Graph Theory}
Let graph $G = (\mathcal{N}, \mathcal{E})$ describe the information exchange among a set $\mathcal{N}$ of agents, where $\mathcal{E} \subset \mathcal{N} \times \mathcal{N}$ is the edge set. If agent $i$ can get information from agent $j$, then $(j, i) \in \mathcal{E}$ and agent $j$ belongs to agent $i$'s neighbour set $\mathcal{N}_{i} = \setc{j}{(j, i)\in \mathcal{E}}$. $G$ is undirected
when $(i, j) \in \mathcal{E}$ if and only if $(j, i) \in \mathcal{E}$. $G$ is connected if there is a path between any two nodes. Let $W = [w_{ij}] \in \reals^{N\times N}$ be the weighted adjacency matrix, with $w_{ij} > 0$ if $j \in \mathcal{N}_{i}$ and $w_{ij} = 0$ otherwise. Let $Deg = \diag(d_{i})_{i\in\mathcal{N}}$, where $d_{i} = \sum_{j=1}^{N} w_{ij}$. Assume that $W = W^{T}$ so the weighted Laplacian of $G$ is $L = Deg - W$. When $G$ is connected and undirected, $0$ is a simple eigenvalue of $L$, $L\mathbf{1}_N = \mathbf{0}_{N}$, $\mathbf{1}_N^{T}L = \mathbf{0}^{T}_{N}$. All other eigenvalues of $L$ are positive, ordered in ascending order as $0 < \lambda_{2}(L) \leq \dots \leq \lambda_{N}(L)$, with $d^* \leq \lambda_N(L) \leq 2 d^*$, where $d^* =\max_i{d_i}$ is the maximal weighted degree of $G$. 
\vspace{-0.25cm}
\section{Game Formulation}\label{sec:formulation}
Consider a group of agents (players) $\mathcal{N} = \{1,\dots, N \}$,
where each player $i \in\mathcal{N}$ controls its local decision (action/strategy) $x_{i} \in \reals^{n}$. Denote $x = col(x_{i})_{i\in\mathcal{N}}\in\reals^{Nn}$ as the decision profile of all agents. Equivalently, we  also write $x = (x_{i}, x_{-i})$ where $x_{-i} = col(\dots, x_{i-1}, x_{i+1}, \dots)$ denotes the decision profile of all agents except player $i$. Agent $i$ aims to optimize
its objective function $\bar{J}_{i}(x_{i},x_{-i})$ (coupled to other players' decisions) with respect to its own decision $x_{i}$ over its
feasible decision set $\Omega_{i}$. Let the globally shared, affine coupled
constrained set be \vspace{-0.25cm}
\begin{align*}
	K \!=\! \prod_{i=1}^{N} \Omega_{i} \cap \setc{x\in\reals^{n}}{\sum_{i=1}^{N}A_{i}x_{i} \leq \sum_{i=1}^{N}b_{i}}
\end{align*}
where $\Omega_{i}\subset \reals^{n}$ is a private feasible set of player $i$, and $A_{i}\in\reals^{m\times n}$, $b_{i}\in\reals^{m}$ is local player information. Let $\Omega = \prod_{i=1}^{N} \Omega_{i}$. 

We focus on average aggregative games where  the cost function $\bar{J}_{i}(x_{i},x_{-i})$ of  each agent $i$ depends on the average of all agents' actions $\sigma(x) \!= \!\frac{1}{N}\sum_{j=1}^{N}x_{j}\!\in \! \reals^{n}$, denoted as $\bar{J}_{i}(x_{i},x_{-i})\!=\!J_{i}(x_{i},\sigma(x))$ to explictly indicate this dependency. In the remainder of the paper we  use either $\bar{J}_{i}$ or $J_{i}$, depending on the context.  Given the other's actions $x_{-i}$, the objective of each player $i$ is to solve the following optimization problem with coupled constraints, \vspace{-0.25cm}
\begin{align}
	\min_{x_{i}}, \, \bar{J}_{i}(x_{i},x_{-i})  \, \, 	s.t. \, \,  (x_{i},x_{-i}) \in K \label{game}
\end{align}
A generalized Nash equilibrium (GNE) is $x^*\!\!=\!col(x^*_i\!)_{i\in  \mathcal{N}}$ such that 
\vspace{-0.35cm}
\begin{align*}
	(\forall i \in \mathcal{N})\qquad x_{i}^{*} {\in} \argmin_{x_{i}} \bar{J}_{i}(x_{i},x_{-i}^{*}), \, \, s.t.\,\, (x_{i},x_{-i}^{*}) \in K \notag
\end{align*}
\begin{asmp} \label{asmp:func}
	For each player $i$, given any  $x_{-i}$, $\bar{J}_{i}(x_{i},x_{-i})$ is continuously differentiable and convex in $x_{i}$ and $\Omega_{i}$ is a compact convex set. The constraint set $K$ is non-empty and satisfies Slater's constraint qualification.
\end{asmp}
Assumption \ref{asmp:func} is a standard assumption which ensures existence of a generalized Nash equilibrium (GNE).

Given the optimization problem (\ref{game}) over $K$ for each agent $i$, let its Lagrangian be defined as 
$
	{\mathcal{L}_{i}}(x_{i},\lambda_{i}; x_{-i}) \!=\! \bar{J}_{i}(x_{i},x_{-i}) + \lambda_{i}^{T}(Ax-b)
$, 
where $A = [A_{1},\dots,A_{N}]$, $b = \sum_{i=1}^{N}b_{i}$ and the dual multiplier is $\lambda_{i}\in \reals^{m}_{+}$. Then, the KKT conditions that an optimal solution $x^{*}_{i}$ with $\lambda_{i}^{*}$ satisfies 
can be written as\vspace{-0.25cm}
\begin{align}\label{kkt_i}
\begin{split}
	\mathbf{0}_{n} &\in \nabla_{x_i}\bar{J}_{i}(x_{i}^{*}, x_{-i}^{*}) + A_{i}^{T}\lambda_{i}^{*} + N_{\Omega_{i}}(x_{i}^{*}) \\
	\mathbf{0}_{m} &\in -(Ax^{*}-b) + N_{\reals^{m}_{+}}(\lambda_{i}^{*}) 
\end{split}
\end{align}
where  
$\nabla_{x_i}\bar{J}_{i}(x_{i},x_{-i}\!) \!=\! \nabla_{x_i}J_{i}(x_{i}, \sigma(x))$, with \vspace{-0.25cm}
$$\! \nabla_{x_i}J_{i}(x_{i}, \sigma(x)\!)\! =\! \nabla_{x_i}J_{i}(x_{i},y\!)|_{y= \sigma(x)} \!+\! \frac{1}{N} \! \nabla_{y}J_{i}(x_{i},y\!)|_{y= \sigma(x)}.$$
and where $N_{\Omega_i}$ is the normal cone operator.

By Theorem 8, \S 4 in \cite{FacchineiGNEP} when $(col(x_{i}^{*})_{i\in\mathcal{N}}, col(\lambda_{i}^{*})_{i\in\mathcal{N}})$ satisfies (\ref{kkt_i}) for all $i\in \mathcal{N}$ then $x^{*}$ is a GNE of the game. Given $x^*$ as a GNE of game, the corresponding Lagrangian multipliers  may be different for the players, i.e.,
$\lambda_1^*\neq \lambda_2^*\neq ,...,\neq\lambda_N^*$. A GNE with the same Lagrangian multiplier for all agents is called {\em variational GNE}, \cite{FacchineiGNEP}, 
with the economic interpretation of no price discrimination and better stability/sensitivity properties,  \cite{shanbhag2}. 
A {\it variational  GNE}  of the game is defined as  $x^* \in K$   solution of the variational inequality $VI(F,K)$:\vspace{-0.25cm}
\begin{align*}
	\inp*{F(x^{*})}{x-x^{*}} \geq 0, \quad \forall x \in K 
\end{align*}
where $F$ denotes  the {\em pseudo-gradient} of the game defined as $F(x) = col(\nabla_{x_i}\bar{J}_{i}(x_{i},x_{-i}))_{i\in\mathcal{N}}\!=\!col(\nabla_{x_i}J_{i}(x_{i}, \sigma(x))_{i\in\mathcal{N}}$, i.e., the stacked vector with the partial gradients for all $i \in \mathcal{N}$. $x^*$ solves $VI(F,K)$ if and only if there  exists a $\lambda^*\in \reals^m$ such that the  KKT conditions are satisfied, \cite[\S 10.1]{FP07}, 
\vspace{-0.25cm}
\begin{equation}
\begin{array}{lll} \label{kkt} 
\mathbf{0}_{Nn}  &\in F(x^*) + A^T \lambda^* +N_{\Omega}(x^*) \\
\mathbf{0}_{m} & \in  -(Ax^*-b) + N_{\reals^m_{+}}(\lambda^*) 
\end{array}
\end{equation}
where $N_{\Omega}(x^*)=\prod_{i=1}^N N_{\Omega_i}(x^*_i)$. Assumption 
\ref{asmp:func}  
 guarantees the existence of a  solution to $VI(F,K)$, by  \cite[Corollary 2.2.5]{FP07}. 
By  \cite[Thm. 9, \S 4]{FacchineiGNEP},  
every solution $x^*$ of  $VI(F,K)$ is a GNE of game. Furthermore, if $x^*$ together with $\lambda^*$ satisfies
 the KKT conditions  (\ref{kkt}) for $VI(F,K)$   then $x^*$  satisfies
the KKT conditions (\ref{kkt_i}) with $\lambda_1^*\!=\!...\!=\!\lambda^*_N\!=\!\lambda^*$,  hence $x^*$ is a variational GNE of game.

Given this aggregative game, our aim is to design a  distributed iterative algorithm that finds a variational GNE under \emph{partial-decision information} over a network. 
\begin{asmp} \label{asmp:pseudo}
$F$ is strongly monotone and Lipschitz continuous, i.e., $\exists \mu \!>\! 0$ and $l_{F} \!>\! 0$ such that $\forall x,\underline{x}$,
$	\inp*{x \!-\!\underline{x}}{F(x)\!-\!F(\underline{x})}\!\geq \!\mu \norm{x\!-\!\underline{x}}^{2}$, and 
$	\norm{F(x)\!-\!F(\underline{x})} \!\leq \! l_{F}\norm{x \!- \!\underline{x}}. 
$
\end{asmp} Assumption 
\ref{asmp:pseudo} is commonly used  in algorithms with fixed-step sizes, 
\cite{NedichDistAgg,FPAuto2016,FPAuto2018,GrammaticoTAC2017,belg,wardrop},  \cite{Shanbhag,parise,PengAgg,Lacra,PengLacraAuto2019},\cite{SalehiShiPavelAuto_2019}, 
and guarantees that a unique variational GNE exists, \cite{FP07}. 
We consider that each agent $i$ does not have information on the other agents actions $x_{-i}$ or on the aggregate value, $\sigma(x)$, and that 
agents communicate only locally with neighbouring agents, over a communication graph $G$. 
\begin{asmp} \label{asmp:graph}
	The communication graph $G=(\mathcal{N},\mathcal{E})$ is undirected and connected.
\end{asmp}
{
\begin{remark}
	In the remainder of the paper, to simplify notation, we assume that $w_{ij} = 1$, $(i,j)\in\mathcal{E}$. 
	The results can be extended immediately to the case $w_{ij} > 0$, $(i,j)\in\mathcal{E}$.
\end{remark}
}
\vspace{-0.25cm}
\section{Distributed Algorithm}\label{sec:alg}

In this section we present our proposed algorithm. 
To offset the lack of full information, each agent $i$ maintains a local estimate $u_i$ of the aggregate $\sigma(x)$ and a local multiplier $\lambda_i$,  and exchanges them with its neighbours over  $G$, in order to learn the true aggregate value and the Lagrange multiplier $\lambda^*$. Each agent also maintains  an additional auxiliary variable $z_i$, used for the coordination of the coupling constraints and to reach consensus of the local multipliers. 

Let  $(x_{i,k}, u_{i,k}, z_{i,k}, \lambda_{i,k})$ denote the tuple with agent $i$'s decision variable $x_{i,k}$, local aggregate estimate $u_{i,k}$, local multiplier $\lambda_{i,k}$, and auxiliary variable $z_{i,k}$ at iteration $k$, respectively.   The goal is that over time each agent will have the same aggregate estimate, equal to the average of the agents actions, the same multiplier, and its decision will correspond to a variational GNE with the coupled constraints met. 
The proposed distributed algorithm is given below. 

\emph{\hspace{-0.4cm}Algorithm 1: }

\noindent\rule{0.49\textwidth}{0.3mm}

$	\hspace{-0.4cm}x_{i,k+1} = P_{\Omega_{i}}[x_{i,k} -\tau_{i}(\nabla_{x_i}J_{i}(x_{i,k},u_{i,k}) + A_{i}^{T}\lambda_{i,k} $

$\hspace{1cm} { + c \sum_{j\in\mathcal{N}_{i}} (u_{i,k}-u_{j,k})} )] $

$	\hspace{-0.4cm}u_{i,k+1} = u_{i,k} - \kappa {c} \sum_{j\in\mathcal{N}_{i}} (u_{i,k}-u_{j,k}) + (x_{i,k+1}-x_{i,k}) $

$	\hspace{-0.4cm}z_{i,k+1} = z_{i,k} + \upsilon_{i} \sum_{j\in\mathcal{N}_{i}}(\lambda_{i,k}-\lambda_{j,k}) \vspace{0.1cm}$

$	\hspace{-0.4cm}\lambda_{i,k+1} \!=\! P_{\reals_{+}^{m}}(\lambda_{i,k} \!-\! \alpha_{i}[\sum_{j\in\mathcal{N}_{i}}(\lambda_{i,k}\!-\!\lambda_{j,k}) \!+\! b_{i}
 -\!A_{i}(2x_{i,k+1}\!$
 
 $\hspace{1cm}-\!x_{i,k}) \!+\! \sum_{j\in\mathcal{N}_{i}}(2z_{i,k+1}\!-\!z_{i,k})\!-\!(2z_{j,k+1}\!-\!z_{j,k})]\! )$
 \noindent\rule{0.49\textwidth}{0.3mm}

\noindent where $x_{i,0}\in\Omega_{i}$, $u_{i,0} = x_{i,0}$, $z_{i,0}, \lambda_{i,0} \in \reals^m$, $\forall i\in \mathcal{N}$, {$c\!>\!0\!$}  and $\tau_{i}, \upsilon_{i},\alpha_{i}, \kappa \in\reals_{++}$ are positive step sizes. 
Note that Algorithm 1 is fully distributed and instead of $\sigma(x)$,  
each agent is using $u_{i}$ (own estimate of the aggregate)  to evaluate  $\nabla_{x_i}J_{i}(x_{i}, u_i) \!=\!    \nabla_{x_i}J_{i}(x_{i},y)|_{y= u_i} \!+\! \frac{1}{N} \! \nabla_{y}J_{i}(x_{i},y)|_{y=u_i}$, its own partial gradient. 

To write the algorithm more compactly, let $x_{k} \!= \!col(\!x_{i,k}\!)_{i\in\mathcal{N}}$, $\!u_{k} \!=\! col(\!u_{i,k}\!)_{i\in\mathcal{N}}$, $\!z_{k} \!=\! col(\!z_{i,k}\!)_{i\in\mathcal{N}}$, $\!\lambda_{k} \!=\! col(\!\lambda_{i,k}\!)_{i\in\mathcal{N}}$. Let $\mathbf{F}(x,u)\!=\! col(\nabla_{x_i}J_{i}(x_{i}, u_{i}))_{i\in\mathcal{N}}$ where $u_{i}\in \reals^{n}$ and $u= col(u_{i})_{i\in\mathcal{N}}$ be the extended pseudo-gradient. Note that when all $u_{i}\!=\!\sigma(x)$,  $\mathbf{F}(x,\mathbf{1}_N \otimes \sigma(x)\!=\! col(\nabla_{x_i}J_{i}(x_{i}, \sigma(x)))_{i\in\mathcal{N}}$, hence   \vspace{-0.22cm} $$\mathbf{F}(x,\mathbf{1}_N \otimes \sigma(x)) \!=\! F(x).$$ Thus, if in $F(x)$ each agent is evaluating the gradient with the true action aggregate value,  in $\mathbf{F}(x,u)$ each agent is using instead its own estimate of the aggregate, $u_{i}$. 

With these we can write Algorithm 1 compactly as, 
\vspace{-0.25cm}\begin{align}
\begin{split}
	x_{k+1} &= P_{\Omega}[x_{k} -\tau(\mathbf{F}(x_{k},u_{k}) + \Lambda^{T}\lambda_{k} {+ c\mathbf{L}_{u}u_{k}} )] \\ \label{eqn:AggIter}
	u_{k+1} &= u_{k} - \kappa { c }\mathbf{L}_{u}u_{k} + (x_{k+1}-x_{k}) \\
	z_{k+1} &= z_{k} + \upsilon \mathbf{L}_{\lambda}\lambda_{k} \\
	\lambda_{k+1} &= P_{\reals_{+}^{Nm}}\left( \lambda_{k} - \alpha[\mathbf{L}_{\lambda}\lambda_{k} + \bar{b} - \Lambda(2x_{k+1}-x_{k}) \right. \\
	&\qquad\qquad\qquad\qquad \left. + \mathbf{L}_{\lambda}(2z_{k+1}-z_{k})] \right)
\end{split}
\end{align}
where $x_{0}\in \Omega \subset \reals^{Nn}$, $u_{0}=x_{0}$, $z_{0} \in  \reals^{Nm}, \lambda_{0} \in\reals^{Nm}_{+}$,  $P_{\Omega} = col(P_{\Omega_{i}})_{i\in\mathcal{N}}$, $P_{\reals_{+}^{Nm}} = col(P_{\reals_{+}^{m}})_{i\in\mathcal{N}}$, $\mathbf{L}_{u} = L \otimes I_{n}$, $\mathbf{L}_{\lambda} = L \otimes I_{m}$, and $L$ is the Laplacian matrix of the graph $G$, $\Lambda = diag(A_{i})_{i\in\mathcal{N}}$, $\bar{b} = col(b_{i})_{i\in\mathcal{N}}$, { $\tau = diag(\tau_{i})_{i\in\mathcal{N}} \otimes I_{n}$, $\upsilon = diag(\upsilon_{i})_{i\in\mathcal{N}}\otimes I_{m}$, and $\alpha = diag(\alpha_{i})_{i\in\mathcal{N}}\otimes I_{m}$}.

\begin{remark} 
{  The algorithm is similar to the algorithm in \cite{Lacra} but it is tailored for aggregative games. } The update for $x_{k}$ is a projected-gradient descent of the local Lagrangian function using player's aggregate estimate $u_{k}$  to evaluate the pseudo-gradient. {Unlike  \cite{Lacra}, here the estimate is separate from the action and needs to track the aggregate of all actions. As a result, 
besides the consensus term it has an additional correction term to account for the action's effect on the average. This acts as a dynamic-tracking term and, as shown in the following, leads to an invariance property of the aggregate estimate, which turns out to be critical for the convergence proof.} The auxiliary variable $z_{k}$ is updated by a discrete-time integration for the multiplier consensual error. Lastly, the update for $\lambda_{k}$ is a combination of a projected-gradient ascent of the local Lagrangian and a proportional-integral term for the multiplier consensus error.
\end{remark}

We prove next the important invariance property that $ \sigma(u_k) = \frac{1}{N}\sum_{i=1}^{N}u_{i,k}$ (the average of all agents' aggregate estimates) is always equal to the actions' true aggregate,  $ \sigma(x_k) = \frac{1}{N}\sum_{i=1}^{N}x_{i,k}$.

\begin{lemma}\label{lemma:Alg3_Prop}
Suppose Assumptions \ref{asmp:func} and \ref{asmp:graph}  hold. Then the following properties hold for iterates $x_{k}, u_{k}, z_{k}, \lambda_k$ generated  by Algorithm 1 or (\ref{eqn:AggIter}). 

(i) $\inp*{\mathbf{1}_N \!\otimes \! I_n \!}{\!u_{k}}\! =\! \inp*{\mathbf{1}_N \! \otimes \! I_n\!}{\!x_{k}}$, hence $\sigma(u_k) \!= \! \sigma(x_k)$, for all $k \geq0$.

(ii) Any {fixed point} $(\overline{x}, \overline{u}, \overline{z}, \overline{\lambda})$ of Algorithm 1 is such $\overline{u}_i=\sigma(\overline{x})$, $\overline{\lambda}_i=\lambda^*$ for all $i \in \mathcal{N}$, where $\overline{x}=x^*$ and $\lambda^*$ satisfy (\ref{kkt}), and $x^*$ is a variational GNE.

\end{lemma}
\emph{ Proof: } See Appendix. \hspace{-0.25cm}

We exploit the invariance in Lemma  \ref{lemma:Alg3_Prop}(i) and, from (\ref{eqn:AggIter}), we next construct an auxiliary iteration with respect to the consensus subspace of the aggregate estimates, which will prove instrumental for the convergence analysis. Let  $C$ denote  the $n$-dimensional consensus subspace for all agents' aggregate estimates, i.e., $	C \!=\! \setc{u}{u = \mathbf{1}_N \otimes c, c\in\reals^{n}}$ and let  $C^\orth$  be its orthogonal complement. 
Note that $C = Null(\mathbf{L}_u) = Range (\mathbf{1}_N \otimes I_n)$ and $C^{\orth} = Range(\mathbf{L}_u)=Null(\mathbf{1}^{T}_N\otimes I_{n})$. 
Any $u\in \reals^{Nn}$ can be decomposed  as $u = P_{\orthc} u + P_{\orth} u$, where $P_{\orthc} u \in C$ and $P_{\orth} u \in C^{\orth}$, by using the projection matrices $P_{\orthc} = \frac{1}{N}\mathbf{1}_N\mathbf{1}^{T}_N\otimes I_{n} = \frac{1}{N} \mathbf{1}_N \otimes \mathbf{1}^{T}_N \otimes I_{n} $ and  $P_{\orth} =I_{Nn}-P_{\orthc}$.  
Note that  $P_{\orthc} u = \mathbf{1}_N \otimes \sigma(u)$,  where  $\sigma(u)=   \frac{1}{N}(\mathbf{1}^{T}_N\otimes I_{n}) u$. 

Any $u_k$  generated by (\ref{eqn:AggIter}) can be decomposed as 
$	u_{k} = P_{\orthc} u_{k} + P_{\orth} u_{k}$ 
where, by using the invariance property in Lemma  \ref{lemma:Alg3_Prop}(i), 
\vspace{-0.25cm}\begin{align}
P_{\orthc} u_{k} = \mathbf{1}_N \otimes \sigma(u_k) = P_{\orthc} x_{k}\quad \forall k\geq 0. \label{eqn:orthEquiv}
\end{align}
Using  this decomposition 
together with (\ref{eqn:AggIter}), consider: 
\vspace{-0.25cm}\begin{align}
\begin{split}
	x_{k+1} &= P_{\Omega}[x_{k} - \tau (\mathbf{F}(x_{k},P_{\orthc}x_{k}+u^{\orth}_{k}) + \Lambda^{T}\lambda_{k}  { + c\mathbf{L}_{u}u^{\orth}_{k}} )] \\ \label{eqn:AggIter2}
	u^{\orth}_{k+1} &= u^{\orth}_{k} - \kappa { c } \mathbf{L}_{u}u^{\orth}_{k} + P_{\orth} (x_{k+1}-x_{k}) \\    
	z_{k+1} &= z_{k} + \upsilon \mathbf{L}_{\lambda}\lambda_{k} \\
	\lambda_{k+1} &= P_{\reals_{+}^{Nm}}\left( \lambda_{k} - \alpha[\mathbf{L}_{\lambda}\lambda_{k} + \bar{b} - \Lambda(2x_{k+1}-x_{k}) \right. \\
	&\qquad\qquad\qquad\qquad \left. + \mathbf{L}_{\lambda}(2z_{k+1}-z_{k})] \right)
\end{split}
\end{align}  where $\!x_{0}\!\in\!\Omega$, $\!u_{0}^{\perp} \!=\!P_{\orth} \!x_{0}$, $\!z_{0} \!\in\!  \reals^{Nm}, \!\lambda_{0} \!\in \!\reals^{Nm}_{+}$. 
The next result relates precisely iterates generated by  Algorithm 1 or (\ref{eqn:AggIter}) to iterates generated by (\ref{eqn:AggIter2}).

\begin{lemma}\label{lemma:Alg3_2_5}
Suppose Assumptions \ref{asmp:func} and \ref{asmp:graph}  hold. Then, any sequence  $\{x_{k}, u_{k}, z_{k}, \lambda_k\}$ generated  by Algorithm 1 or (\ref{eqn:AggIter}) with initial conditions $x_0, u_0\!=\!x_0, z_0, \lambda_0$ can be  derived from some sequence $\{x'_{k}, u^{'\orth}_{k}, z'_k, \lambda'_k\}$ generated by (\ref{eqn:AggIter2}) with initial conditions $x'_0\!=\!x_0, u^{'\orth}_0\!=\!P_{\orth} x_0, z'_0\!=\!z_0, \lambda'_0\!=\!\lambda_0$, as in \vspace{-0.25cm}
\begin{align}
x_{k}\! =\! x'_{k}, u_k\! =\! P_{\orthc}x'_{k}\!+\!u^{'\orth}_{k}, z_k\!=\!z'_{k},\lambda_k\!=\! \lambda'_k. \label{eqn:Alg5_to_3}
\end{align}
\end{lemma}

\emph{ Proof: }
The proof follows an induction argument. 
Due to the  initial conditions, $P_{\orthc}x'_{0} \! +\!u^{'\orth}_{0} \!=P_{\orthc}x_{0} \! +\! \!P_{\orth} x_0 \! =\!x_0 \!=\!u_0$, hence (\ref{eqn:Alg5_to_3}) holds for $k=0$. 

Suppose  (\ref{eqn:Alg5_to_3}) holds at step $k$. Then,   from $ z'_{k+1}\!=\!z'_{k} \!+\! \upsilon \mathbf{L}_{\lambda}\lambda'_{k} $ (cf.   (\ref{eqn:AggIter2})) with $z'_k\!=\!z_{k},\lambda'_k\!=\! \lambda_k$, it follows that $ z'_{k+1}\!= \!z_{k} \!+\! \upsilon \mathbf{L}_{\lambda}\lambda_{k}$ hence by the $z$-update in  (\ref{eqn:AggIter}), $z'_{k+1}\! =\! z_{k+1}$. Next, using (\ref{eqn:Alg5_to_3}) into  the right-hand side of the $x'$-update in  (\ref{eqn:AggIter2}), yields  cf. (\ref{eqn:AggIter}),\vspace{-0.25cm}
\begin{align}\label{eq:xprime_kp1}
x'_{k+1} \!&=\! P_{\Omega}[x'_{k} \!-\! \tau (\mathbf{F}(x'_{k},P_{\orthc}x'_{k}\!+\!u^{'\orth}_{k}) 
\nonumber \\
&\qquad + \Lambda^{T}\lambda'_{k} { + c\mathbf{L}_{u}u^{'\orth}_{k}} )] \\
&=\! P_{\Omega}[x_{k} \!-\! \tau(\mathbf{F}(x_{k},\! u_k) \!+\! \Lambda^{T}\lambda_{k} { + c\mathbf{L}_{u}u_{k}} )] \nonumber \\
 &= x_{k+1}. \nonumber
\end{align}Hence the first and third relations in  (\ref{eqn:Alg5_to_3}) hold at step $k+1$, and using them on the right-hand side of the $\lambda'$-update in (\ref{eqn:AggIter2}) yields $\lambda'_{k+1}\!=\! \lambda_{k+1}$, where $\lambda_{k+1}$ is generated by (\ref{eqn:AggIter}).

Lastly, we  show the second relation in  (\!\ref{eqn:Alg5_to_3}) at step $k\!+\!1\!$. Thus, \vspace{-0.2cm}
\begin{align*}
&P_{\orthc}x'_{k+1}\!+\!u^{'\orth}_{k+1}\!=\!P_{\orthc}x_{k+1}\! +\! u^{'\orth}_{k} \!-\! \kappa {c} \mathbf{L}_{u}u^{'\orth}_{k} \!+\! P_{\orth} (x'_{k+1}\!-\!x'_{k}) \\
&\!=\! P_{\orthc}x_{k+1}\! +\! u_k\! -\! P_{\orthc}x_{k}\!-\! \kappa {c} \mathbf{L}_{u} (u_k\! -\! P_{\orthc}x_{k}) \!+\! P_{\orth} (x_{k+1}\!-\!x_{k}),
\end{align*}where we used relations for $x$ at step $k+1$ cf. (\ref{eq:xprime_kp1}) and for $u^{'\orth}$ at step $k$ cf.  (\ref{eqn:Alg5_to_3}). Using $P_{\orthc}x_{k} \!+\!  P_{\orth} x_{k} \!=\!  x_{k}$, $\mathbf{L}_{u} P_{\orthc}x_{k}\!=\!\mathbf{0}$ it follows that $P_{\orthc}x'_{k+1}\!+\!u^{'\orth}_{k+1}\!=\! \! u_k\! -\! \kappa {c} \mathbf{L}_{u}u_{k} \!+\!  (x_{k+1}\!-\!x_{k})\!=\! u_{k+1}$, cf. the $u$-update in (\ref{eqn:AggIter}), 
and the argument is complete. 
\hfill $\Box$ 
Lemma \ref{lemma:Alg3_2_5} is instrumental in what follows. Based on it, convergence of Algorithm 1 or (\ref{eqn:AggIter}) can be established once convergence of  (\ref{eqn:AggIter2}) is shown. Note that any $u_{k}^{\orth}$ generated by (\ref{eqn:AggIter2}) is such that $u_{k}^{\orth} \in C^\orth$, for all $k \geq0$. 
 We show next that  (\ref{eqn:AggIter2}) can be written as a preconditioned forward-backward iteration, \vspace{-0.25cm}
\begin{equation}
	\mathbf{0}_{2N(n+m)}\in \mathcal{A}\varpi_{k} + \mathcal{B}\varpi_{k+1} + \Phi(\varpi_{k+1} - \varpi_{k}) \label{eqn:dup}
\end{equation}
where $\varpi \!= \! (x,u^{\orth},z,\lambda) \in \reals^{Nn} \! \times \!C^\orth  \!\times  \!\reals^{Nm} \!\times \! \reals^{Nm}_+$, for two operators $\mathcal{A}$ and $ \mathcal{B}$ and  matrix $\Phi$ defined as, \vspace{-0.25cm} 
\begin{align}
\begin{split} \label{eqn:equDynOrth}
\mathcal{A} &: \varpi  \!\mapsto  \! \begin{bmatrix}
\mathbf{F}(x,P_{\orthc}x \!+ \!u^{\orth}) \\
{c}\mathbf{L}_{u}u^{\orth} \\
\mathbf{0}_{Nm} \\
\mathbf{L}_{\lambda} \lambda + \bar{b}
\! \end{bmatrix} \\
\mathcal{B} &: \varpi \! \mapsto \! \begin{bmatrix} \!
\!N_{\Omega}(x)\! \\
\!\mathbf{0}_{Nn}\! \\
\!\mathbf{0}_{Nm}\! \\
\!N_{\reals_{+}^{Nm}(\lambda) \!}
\! \end{bmatrix}  + \begin{bmatrix} \!
\mathbf{0} &\!\! { \mathbf{0}} & \!\mathbf{0} &\! \Lambda^T \\
{ \mathbf{0}} &\! \!\mathbf{0} &\! \mathbf{0} &\! \mathbf{0} \\
\mathbf{0} & \!\!\mathbf{0} & \!\mathbf{0} &\! -\mathbf{L}_{\lambda} \\
-\Lambda &\! \!\mathbf{0} & \!\mathbf{L}_{\lambda} &\! \mathbf{0}
\!\end{bmatrix} \varpi \\
\Phi &= \begin{bmatrix}
	\tau^{-1}   { + \kappa^{-1}P_{\orth} } & { -\kappa^{-1}P_{\orth} } & \mathbf{0} & -\Lambda^T \\
	{-\kappa^{-1}P_{\orth}} & \kappa^{-1}{I_{Nn}} & \mathbf{0} & \mathbf{0} \\
	\mathbf{0} & \mathbf{0} & \upsilon^{-1} & \mathbf{L}_{\lambda} \\
	-\Lambda & \mathbf{0} & \mathbf{L}_{\lambda} & \alpha^{-1}	 
	\end{bmatrix}
\end{split}
\end{align}
{ where  $\mathbf{0}$ is an all-zero matrix of appropriate dimension.}
\begin{lemma} \label{lemma:equivRep}Suppose that Assumptions \ref{asmp:func}, \ref{asmp:pseudo} and \ref{asmp:graph} hold and 
let $\varpi_{k} = col(x_{k},u_{k}^{\orth},z_{k},\lambda_{k})$, $\mathcal{A}$, $\mathcal{B}$ and $\Phi$ be defined as in (\ref{eqn:equDynOrth}). Suppose that $\Phi \succ 0$ and $\Phi^{-1}\mathcal{B}$ is maximally monotone. Then the following hold:

\noindent (i): Iterates (\ref{eqn:AggIter2}) are equivalently written as (\ref{eqn:dup}), i.e., \vspace{-0.1cm} \begin{equation} \label{eqn:optFBequiv}
\varpi_{k+1} = (\Id + \Phi^{-1}\mathcal{B})^{-1}\circ(\Id - \Phi^{-1}\mathcal{A})\varpi_{k} = T_{2}\circ T_{1}\varpi_{k}
\end{equation}
where $T_{2} = (\Id + \Phi^{-1}\mathcal{B})^{-1}$ and $T_{1} = \Id - \Phi^{-1}\mathcal{A}$.

\noindent (ii) Any {fixed point} of (\ref{eqn:AggIter2}) is a zero  of $\mathcal{A} \!+\! \mathcal{B}$ and a fixed point of $T_{2} \!\circ \! T_{1}$. Furthermore, any such point  $\bar{\varpi} \! =\! col(\bar{x},\bar{u}^{\orth},\bar{z},\bar{\lambda})$ satisfies $\bar{x}\!=\!x^*$, $\bar{u}^{\orth}\!=\!\mathbf{0}_{Nn}$, $\bar{\lambda} \!=\! \mathbf{1}_N \! \!\otimes \! \!\lambda^*$,  where $x^*$ and  $\lambda^*$  satisfy the KKT conditions (\ref{kkt}), hence $x^*$ is a variational GNE.
\end{lemma}
\emph{ Proof: } 
(i)  The equivalence {for $u_{k+1}^{\orth}$, $z_{k+1}$, and $\lambda_{k+1}$}  between  (\ref{eqn:dup})  and (\ref{eqn:AggIter2}) can be shown by expanding  (\ref{eqn:dup}) with $\mathcal{A}$, $\mathcal{B}$ and $\Phi$  as in (\ref{eqn:equDynOrth}), cancelling terms.  { Similarly, for the $x_{k+1}$ update,   by expanding (\ref{eqn:dup}), replacing the $u_{k+1}^{\orth}$ term with $u^{\orth}_{k} - \kappa c \mathbf{L}_{u}u^{\orth}_{k} + P_{\orth} (x_{k+1}-x_{k})$ cf. (\ref{eqn:AggIter2}), using $P_{\orth} P_{\orth} \!=\!P_{\orth} $,}  cancelling terms  and using  $\tau^{-1} N_{\Omega}(x)=N_\Omega(x)$, $P_\Omega \!=\! (Id + N_\Omega)^{-1}$ yields (\ref{eqn:AggIter2}). 
Since $\Phi  \!\succ \!0$ by assumption, (\ref{eqn:dup}) is equivalent to $	(\Id \!-\! \Phi^{-1}\mathcal{A})(\varpi_{k}) \!\in \! (\Id  \!+\! \Phi^{-1}\mathcal{B})(\varpi_{k+1})$, which  can be written as (\ref{eqn:optFBequiv}), based on   the fact that  $(\Id + \Phi^{-1}\mathcal{B})^{-1}$ is singled-valued  (cf. \cite[Prop. 23.7]{monoBookv1} by $\Phi^{-1}\mathcal{B}$ maximally monotone).

(ii): Let $\bar{\varpi} = col(\bar{x},\bar{u}^{\orth},\bar{z},\bar{\lambda})$ be the {fixed point} of (\ref{eqn:AggIter2}) or (\ref{eqn:dup}), hence is a fixed point of $T_{2}\circ T_{1}$. By continuity,  the following equivalences hold, 
$\bar{\varpi} = (\Id + \Phi^{-1}\mathcal{B})^{-1}\circ(\Id - \Phi^{-1}\mathcal{A})\bar{\varpi} 
\Leftrightarrow(\Id + \Phi^{-1}\mathcal{B})(\bar{\varpi}) \in (\Id - \Phi^{-1}\mathcal{A})(\bar{\varpi}) 
\Leftrightarrow \mathbf{0}_{2N(n+m)}\in (\mathcal{B} + \mathcal{A})(\bar{\varpi})$. Thus, using  (\ref{eqn:equDynOrth}), it follows that $\bar{\varpi} = col(\bar{x},\bar{u}^{\orth},\bar{z},\bar{\lambda})$ satisfies: 
$\mathbf{0}_{Nn} \!\in \!  N_{\Omega}(\bar{x}) \!+\! \mathbf{F}(\bar{x},P_{\orthc}\bar{x} \!+ \! \bar{u}^{\orth}) + \Lambda^{T}\bar{\lambda}$, 
$\mathbf{0}_{Nn} \!=\!\mathbf{L}_{u} \bar{u}^{\orth}$, 
$\mathbf{0}_{Nm} \!= \!\mathbf{L}_{\lambda} \bar{\lambda}$, and  $\mathbf{0}_{Nm} \! \in \!N_{\reals_{+}^{Nm}}(\bar{\lambda}) \!+\! \mathbf{L}_{\lambda} \bar{\lambda} \! +\! \bar{b} \!-\! \Lambda \bar{x} \!+\! \mathbf{L}_{\lambda} \bar{z}. $
By Assumption \ref{asmp:graph}, the second and third relations imply that $\bar{u}^{\orth}=\mathbf{0}_{Nm}$ (since  $\bar{u}^{\orth} \in C^\orth$) and $ \bar{\lambda} = \mathbf{1}_N \otimes \lambda^*$, for a $\lambda^* \in \reals^m$.  Using these in the first and fourth relations together with { $P_{\orthc}\bar{x} = \mathbf{1}_{N}\otimes \sigma(\bar{x})$},  leads to  $
\mathbf{0}_{Nn} \in  N_{\Omega}(\bar{x}) + \mathbf{F}(\bar{x},\mathbf{1}_N \otimes \sigma(\bar{x})) + A^{T}\lambda^* 
$
which, using $\mathbf{F}(\bar{x},\mathbf{1}_N \otimes \sigma(\bar{x})=F(\bar{x})$, is the first relation in (\ref{kkt}). From 
$\mathbf{0}_{Nm} \! \in  \!N_{\reals_{+}^{Nm}}(\mathbf{1}_N \otimes \lambda^*) \!+\! \bar{b} \!-\! \Lambda \bar{x} \!+\! \mathbf{L}_{\lambda} \bar{z}$,  premultiplying by $ (\mathbf{1}^T \!\otimes \! I_m)$
as in the proof of Lemma  \ref{lemma:Alg3_Prop}(ii), the second one in (\ref{kkt})holds, so  $\bar{x}=x^*$ is a variational GNE.\hfill $\Box$

\begin{remark} Lemma \ref{lemma:Alg3_2_5} and \ref{lemma:equivRep} show that Algorithm 1 is related to a forward-backward iteration (\ref{eqn:dup}) or (\ref{eqn:optFBequiv}), \cite{monoBookv1}, with a preconditioning matrix $\Phi$. 
The symmetric matrix $\Phi$ is needed, as pointed out in \cite{PengLacraAuto2019} and \cite{Lacra}, to be able to distribute the backward step. Our proof techniques are similar  to those in \cite{Lacra}, however,  the operators $\mathcal{A}$  and $\mathcal{B}$  are slightly different, while the metric matrix  $\Phi$ has a different structure, with a new block involving the projection matrix $P_\orth$ which needs handled separately. One of the reasons is that in the general setup in  \cite{Lacra} the action is part of the estimate vector, while in our setup  the actions and the aggregate estimate are separate and only consensus on the aggregate is needed. 
\end{remark}
\vspace{-0.25cm}
\section{Convergence Analysis}\label{sec:convergence}

In this section we show convergence of Algorithm 1. Based on Lemma 
\ref{lemma:Alg3_2_5} its convergence can be established once convergence of  (\ref{eqn:AggIter2}) is shown. In turn,  (\ref{eqn:AggIter2})  is equivalent to  (\ref{eqn:dup}) or (\ref{eqn:optFBequiv}), (cf. Lemma  \ref{lemma:equivRep}).
Convergence of  iteration (\ref{eqn:optFBequiv})  is guaranteed when $\Phi^{-1}\mathcal{A}$ is a cocoercive and $\Phi^{-1}\mathcal{B}$ is a monotone operator, cf. Theorem 25.8 in \cite{monoBookv1}.
However,  $\mathcal{A}$ is defined in terms of the extended pseudo-gradient  $\mathbf{F}$ for which monotonicity properties are not guaranteed to hold on the augmented space of actions and aggregate estimates (only strong monotonicity of $F(x)$ is assumed).  Instead, we will show that $\mathcal{A}$ satisfies a \emph{restrictive cocoercive } property and  that $\mathcal{B}$ is maximally monotone, and  then similar properties for $\Phi^{-1}\mathcal{A}$ and $\Phi^{-1}\mathcal{B}$. This turns out to be sufficient to prove convergence because the restrictive property is with respect to the aggregate consensus subspace where the zeros  of $\mathcal{A} + \mathcal{B}$ lie, cf. Lemma  \ref{lemma:equivRep}(ii). 
 To show the  \emph{restrictive cocoercive } property of $\mathcal{A}$ we balance the strong monotonicity of the pseudo-gradient $F(x)$ (on the aggregate consensus subspace $C$) with that of the Laplacian on its orthogonal component $C^\orth$, under a Lipschitz assumption on the extended pseudo-gradient $\mathbf{F}(x,u)$.
 \begin{asmp}\label{asmp:lip}
	\hspace{-0.2cm}  The extended pseudo-gradient $\!\mathbf{F}\!(x,\!u\!)$ is Lipschitz continuous in both arguments, i.e., $\!\exists l_{F}^{u}\!\!>\!\!0$, $\!l_{F}^{x}\!\!>\!\!0$ s.t. 

$	\norm{\mathbf{F}(x,u)-\mathbf{F}(x,\underline{u})} \leq l_{F}^{u}\norm{u-\underline{u}} \quad \forall x, u, \underline{u} \in \reals^{Nn}, 
$ 

$	\norm{\mathbf{F}(x,u)-\mathbf{F}(\underline{x},u)} \leq l_{F}^{x}\norm{x-\underline{x}} \quad \forall x, \underline{x}, u \in \reals^{Nn}.
$
\end{asmp} Assumption \ref{asmp:lip} is not restrictive and is also used in other distributed (G)NE algorithms over networks, e.g. \cite{NedichDistAgg}, \cite{parise},  \cite{SalehiShiPavelAuto_2019}. { A sufficient condition for it to hold is that $\mathbf{F}$ is C$^1$ with bounded Jacobian. For quadratic games, it is automatically satisfied.  }

The following lemma establishes a (restricted) monotonicity property on part of the $\mathcal{A}$ operator, (\ref{eqn:equDynOrth}), denoted $\tilde{\mathcal{A}}$. 

\begin{lemma} \label{lemma:restrictCoer}
Let $\tilde{\mathcal{A}}$ be defined as, \vspace{-0.22cm}
\begin{align}
	\tilde{\mathcal{A}} &: (x,u^{\orth})
\mapsto \begin{bmatrix}
	\mathbf{F}(x,P_{\orthc}x + u^{\orth}) \\
	{c}\mathbf{L}_{u} u^{\orth}
	\end{bmatrix} \label{tildeAdef}
\end{align}
where $(x,u^{\orth}) \! \in \! \reals^{Nn} \! \times \!C^\orth$. Suppose that Assumptions \ref{asmp:func} to 
\ref{asmp:lip} hold. Then, for any $(x,u^{\orth})$ 
and any $(\underline{x},\underline{u}^{\orth})$ with $\underline{u}^{\orth}=\mathbf{0}_{Nn}$,   
\vspace{-0.25cm}
\begin{align}
	\inp*{\!\begin{bmatrix}\!
	x-\underline{x} \\
	u^{\orth} - \underline{u}^{\orth}\!
\end{bmatrix}\!}{\!\tilde{\mathcal{A}}(x,u^{\orth}) \!-\! \tilde{\mathcal{A}}(\underline{x}, \underline{u}^{\orth})\!} \!\geq \! \mu_{\tilde{\mathcal{A}}} \norm{\!\begin{matrix}
	x-\underline{x} \\
	u^{\orth} - \underline{u}^{\orth}
\!\end{matrix}}^{2} \label{eq:resUmono}
\end{align}
where \vspace{-0.35cm}
\begin{equation*}
	\mu_{\tilde{\mathcal{A}}} = \lambda_{\min}\left( \begin{bmatrix}
	\mu & \frac{-l^{u}_{{F}}}{2}  \\
	\frac{-l^{u}_{{F}}}{2} & {c}\lambda_{2}(L)
	\end{bmatrix}\right).
\end{equation*}
Furthermore, $\tilde{\mathcal{A}}$ is restricted monotone if ${c}\lambda_{2}(L) \geq \frac{(l^{u}_{F})^{2}}{4\mu}$ and strongly monotone if the inequality is strict.
\end{lemma}
\emph{ Proof: } 
With $\tilde{\mathcal{A}}$ as in (\ref{tildeAdef}),  for any $(x,\!u^{\orth}\!) $, 
$(\underline{x},\underline{u}^{\orth}\!)$ with $\!\underline{u}^{\orth}\!=\!\mathbf{0}_{Nn}$,  
the left-hand side of  (\ref{eq:resUmono}) is  written as \vspace{-0.22cm}
\begin{align}\label{eq:atilde_1}
	&\inp*{\!\begin{bmatrix}\!
	x-\underline{x} \\
	u^{\orth} - \underline{u}^{\orth}\!
\end{bmatrix}\!}{\!\tilde{\mathcal{A}}(x,u^{\orth}) \!-\! \tilde{\mathcal{A}}(\underline{x}, \underline{u}^{\orth})\!} 
\\
	& 
=	\inp*{\!x \!-\!\underline{x}\!}{\!\!\mathbf{F}(x,\!P_{\orthc}x \!+\! u^{\orth}\!)\!} 
\!-\! \inp*{\!x \!-\!\underline{x}\!}{\!\!\mathbf{F}(\underline{x},\!P_{\orthc}\underline{x}\!)\!} \nonumber
\\
	& 
	\!+\! \inp*{\!u^{\orth} \!\!-\! \underline{u}^{\orth}\!\!}{\!\!{c}\mathbf{L}_{u}\!(\!u^{\orth} \!-\! \underline{u}^{\orth}\!)\!} \nonumber
	\!+\! \inp*{\!x \!-\!\underline{x}\!}{\!\!\mathbf{F}(\!x,\!P_{\orthc}x\!)\!}\! 
\!-\! \inp*{\!x \!-\!\underline{x}\!}{\!\!\mathbf{F}(x,\!P_{\orthc}x\!)\!}.  \nonumber
\end{align}
Note that 
$\inp*{\!x \!-\!\underline{x}\!}{\!\!\mathbf{F}\!(x,\!P_{\orthc}x) \!-\! \mathbf{F}\!(\underline{x},\!P_{\orthc}\underline{x}\!)\!} 
\!	\geq
\!	 \mu \norm{\!x\! -\!\underline{x}\!}^{\!2}, 
$ by $\mathbf{F}(x,\!P_{\orthc}x)\! =\! \mathbf{F} (x,\!\mathbf{1}_N \!\otimes \!\sigma(x)) \! =\! F(x)$ and Assumption \ref{asmp:pseudo}. Also,   \vspace{-0.22cm}
$$\hspace{-0.2cm}\inp*{\!x\! -\!\underline{x}\!}{\!\!\mathbf{F}(x,\!P_{\orthc}x \!+\! u^{\orth}\!) \!-\! \mathbf{F}(x,\!P_{\orthc}x\!)\!} 
	\!\geq \!-l^{u}_{\mathbf{F}}\norm{u^{\orth}  \!-\! \underline{u}^{\orth}\!}\!\!\norm{x \!-\!\underline{x}}, 
$$ by Assumption \ref{asmp:lip}. Since $u^\orth \!\in \!C^\orth$, $\underline{u}^{\orth}\!=\!\mathbf{0}_{Nn}$,  \vspace{-0.22cm}
$$
	\inp*{u^{\orth} \!-\! \underline{u}^{\orth}\!}{{c}\mathbf{L}_{u}(u^{\orth} \!-\! \underline{u}^{\orth}\!)\!} 
	\!\geq {c}\lambda_{2}(L)\norm{u^{\orth}\! -\! \underline{u}^{\orth}\!}^{2},
$$ by  Assumption  \ref{asmp:graph}. Using these inequalities  
into (\ref{eq:atilde_1}) yields  \vspace{-0.22cm}
\begin{align*}
	&\inp*{\!\begin{bmatrix}\!
	x-\underline{x} \\
	u^{\orth} - \underline{u}^{\orth}\!
\end{bmatrix}\!}{\!\tilde{\mathcal{A}}(x,u^{\orth}) \!-\! \tilde{\mathcal{A}}(\underline{x}, \underline{u}^{\orth})\!} \\
	& \quad \geq	 \begin{bmatrix}
	\norm{x-\underline{x}} & \!\!\! \norm{u^{\orth}-\underline{u}^{\orth}}
	\end{bmatrix}
	\begin{bmatrix}
	\mu & \frac{-l^{u}_{F}}{2}  \\
	\frac{-l^{u}_{F}}{2} & \! {c}\lambda_{2}(L)\!
	\end{bmatrix} \begin{bmatrix}
		\norm{x-\underline{x}} \\
		\norm{u^{\orth}-\underline{u}^{\orth}}\!
	\end{bmatrix}
\end{align*}
from which (\ref{eq:resUmono}) follows. 
\hfill $\Box$
\begin{remark}{The condition ${c}\lambda_2(L)\!>\!\frac{(l_F^u)^2}{4\mu}$ in Lemma \ref{lemma:restrictCoer} is critical to ensure the
restricted monotonicity of $\tilde{\mathcal{A}}$ and  the convergence of the proposed method, and the parameter ${c}$ provides the extra degree of freedom to satisfy it. }
Lemma \ref{lemma:restrictCoer} shows that $\tilde{\mathcal{A}}$ is monotone when restricted to a subspace where $\underline{u}^{\orth} =\mathbf{0}_{Nn}$, i.e.,  when  $\underline{u} = \underline{u}^{\orthc}\in C$. This means that restricted monotonicity  is with respect to the $n$-dimensional aggregate consensus subspace (independent of $N$). This is unlike  \cite{Lacra} where the restricted monotonicity  is with respect to the $Nn$-dimensional consensus subspace (cf. Lemma 3 in \cite{Lacra}). In fact, the bound in  (\ref{eq:resUmono}) is tighter than the one in Lemma 3 in \cite{Lacra} (or Theorem 5 in \cite{Gadjov})  and does not depend on $N$ (number of agents). {We emphasize that we are able to obtain this better bound due to the invariance property of the aggregate estimate. The invariance property is induced by the extra term in the aggregate estimate update in Algorithm 1, that is not present in the algorithm in \cite{Lacra}. }
\end{remark}

Using the fact that $\tilde{\mathcal{A}}$ is restricted monotone we can now show properties for  the operators $\mathcal{A}$ and $\mathcal{B}$, (\ref{eqn:equDynOrth}).

\begin{lemma} \label{lemma:restrictCo}
Suppose Assumptions \ref{asmp:func} 
to \ref{asmp:lip} hold and $\!{c}\lambda_{2}\!(L) \!>\frac{\! (l^{u}_{F}\!)^{\!2}}{4\mu}\!\!$.
	Then the following hold for operators $\mathcal{A}$ and $\mathcal{B}$, (\ref{eqn:equDynOrth}).

\noindent (i) $\mathcal{B}$ is maximally monotone.

\noindent (ii) $\mathcal{A}$ is $\beta$-restricted cocoercive: for any $\varpi \!=\! (x,u^{\orth}\!,z,\lambda)$  with $u^{\orth} \!\in \!C^\orth$ and any $\underline{\varpi} \!=\! (\underline{x},\underline{u}^{\orth}\!,\underline{z},\underline{\lambda})$ with $\underline{u}^{\orth} \!=\! \mathbf{0}_{Nn}$, the following holds \vspace{-0.25cm}
\begin{align*}
	\inp*{\varpi - \underline{\varpi}}{\mathcal{A}\varpi - \mathcal{A}\underline{\varpi}} \geq \beta\norm{\mathcal{A}\varpi - \mathcal{A}\underline{\varpi}}^{2}
\end{align*}
where $\beta \! \!=\! \min\set{\frac{\mu_{\tilde{\mathcal{A}}}}{\theta^2}, \!\frac{1}{\lambda_{N}(L)}}$, $\!\theta^2\! \!=\! \max\set{(l_{F}^{x})^{2},\! (l_{F}^{u})^{2}\!+\!({c}\lambda_{N}(L))^{2}}$.
\end{lemma}
\emph{ Proof: } 
	(i):  $\mathcal{B}$ is written as the sum of two operators, one being a Cartesian product of normal cone operators, hence maximally monotone operator, and the other one 
being a skew-symmetric matrix, hence also maximally monotone, with full domain. 
Thus  $\mathcal{B}$ itself is maximally monotone, \cite{monoBookv1}.

(ii):  Let $\varpi \!= \!(x,u^{\orth} ,z,\lambda)$, $\underline{\varpi} \!=\! (\underline{x}, \underline{u}^{\orth},\underline{z}, \underline{\lambda})$ and $\tilde{\mathcal{A}}$   (\ref{tildeAdef}). Then  
\vspace{-0.22cm}
\begin{align}\label{eq:A_prop1}
	&\inp*{\!\varpi \!-\! \underline{\varpi}\!}{\!\mathcal{A}\varpi \!-\! \mathcal{A}\underline{\varpi}\!} 
\! =\!	\inp*{\!\begin{bmatrix}\!
	x-\underline{x} \\
	u^{\orth} - \underline{u}^{\orth}\!
\end{bmatrix}\!}{\!\tilde{\mathcal{A}}(x,u^{\orth}) \!-\! \tilde{\mathcal{A}}(\underline{x}, \underline{u}^{\orth})\!} \nonumber \\
	& \qquad + \inp*{\lambda - \underline{\lambda}}{\mathbf{L}_{\lambda}(\lambda - \underline{\lambda})}   
	\end{align}
Using (\ref{eq:resUmono}) 
and $	\inp*{\lambda \!-\! \underline{\lambda}}{\mathbf{L}_{\lambda}(\lambda \!-\! \underline{\lambda})} \!\geq \!\frac{1}{{\lambda_{N}(L)}}\norm{\mathbf{L}_{\lambda}(\lambda \!-\! \underline{\lambda})}^{2} 
$ for the second term,  it follows that for any $u^{\orth}\! \in\! C^\orth$ and $\underline{u}^{\orth} \!=\! \mathbf{0}_{Nn}$,   \vspace{-0.2cm}
\begin{align*}\label{eq:A_prop2}
	&\inp*{\varpi \!-\! \underline{\varpi}}{\mathcal{A}\varpi \!-\! \mathcal{A}\underline{\varpi}} 
\! \geq \! \mu_{\tilde{\mathcal{A}}} \norm{\!
	\begin{bmatrix}
	x -\underline{x} \\
	u^{\orth} \!-\! \underline{u}^{\orth}
	\!\end{bmatrix}}^2 
\!+\! \frac{1}{\lambda_{N}(L)}\norm{\mathbf{L}_{\lambda}(\lambda \!-\! \underline{\lambda})}^{2}.
	\end{align*}
On the other hand, using $\tilde{\mathcal{A}}$ as in  (\ref{tildeAdef}), Lipschitz properties of $\mathbf{F}$, bounds on the eigenvalues of the Laplacian matrix $L$,  we obtain 
 \vspace{-0.22cm}
\begin{align*}
	&\norm{ 
\tilde{\mathcal{A}}(x,u^{\orth}) - \tilde{\mathcal{A}}(\underline{x}, \underline{u}^{\orth})}^2 
	\leq \theta^2
	 \norm{
	\begin{bmatrix}
	x -\underline{x} \\
	u^{\orth} - \underline{u}^{\orth}
	\end{bmatrix}}^2 
\end{align*}
where  
 $\theta^2 \!\!=\! \max\set{(l_{F}^{x}\!)^{2}\!, \!(l_{F}^{u}\!)^{2}\!+\! ({c}\lambda_{N}(L)\!)^{2}}$. Using this in the foregoing 
leads to the inequality  in (ii). \hfill $\Box$

The following lemma shows how agents can select step sizes independently such that $\Phi \succ 0$. 
\begin{lemma} \label{lemma:PhiEig}
Given any $\delta \!>\! 0$ and $\kappa<\frac{1}{\delta}$, if step sizes are selected such that \vspace{-0.2cm}
\begin{align*}
	\tau_{i} &\leq \frac{1}{
	\max_{j=1,\dots,n}\set{\sum^{m}_{k=1}|[A^{T}_{i}]_{jk}|} + \delta + {\frac{1}{\kappa(1-\kappa\delta)}}} \\ 
	\upsilon_{i} &\leq  (2d_{i}+\delta)^{-1} \\ \vspace{-0.25cm} 
	\alpha_{i} &\leq (\max_{j=1,\dots,m}\set{\sum^{n}_{k=1}|[A_{i}]_{jk}|} + 2d_{i} + \delta)^{-1} 
	\end{align*}
where $d_{i} = |\mathcal{N}_{i}|$, 
 then the matrix  $\Phi - \delta I \succeq 0$ and $\Phi \succ 0$. 
\end{lemma}
{\emph{ Proof: } 	See the Appendix. \hfill $\Box$}

The next result shows that $\Phi^{-1}\mathcal{B}$ and $\Phi^{-1}\mathcal{A}$ satisfy a monotonicity property in the $\Phi$-induced norm.

\begin{lemma} \label{lemma:PhiNormMono}
 Suppose Assumptions \ref{asmp:func} 
to \ref{asmp:lip} hold and ${c}\lambda_{2}\!(L) \!>\frac{\! (l^{u}_{F})^{\!2}}{4\mu}\!\!$. Take any $\delta > \frac{1}{2\beta}$ where $\beta$ is from Lemma \ref{lemma:restrictCo}, $\kappa \!<\!\frac{1}{\delta}$ and  step sizes $\tau_{i}$, $\upsilon_{i}$ and $\alpha_{i}$ chosen to satisfy Lemma  \ref{lemma:PhiEig}. 
Then under the $\Phi$-induced norm $\norm{\cdot}_{\Phi}$ the following hold:

\noindent (i)	 $\!\Phi^{-1} \!\mathcal{B}$ is maximally monotone and $T_{2} \!= \!(\!I\! +\! \Phi^{-1}\!\mathcal{B})^{\!-1}\! \!\in \! \mathfrak{A}(\!\frac{1}{2}\!)$.

\noindent (ii)	$\!\Phi^{-1}\!\mathcal{A}$ is $\beta\delta$-restricted cocoercive and  $T_{1} = Id-\Phi^{-1}\mathcal{A}$  is restricted nonexpansive and furthermore, the following holds  for any $\varpi = (x,u^{\orth},z,\lambda)$ and any $\underline{\varpi} = (\underline{x},\mathbf{0}_{Nn},\underline{z},\underline{\lambda})$, \vspace{-0.22cm}
\begin{align*}
\!	\norm{\!T_{1}\varpi \!-\! T_{1}\underline{\varpi}\!}_{\Phi}^{2} \!\leq &\!\norm{\varpi \!-\! \underline{\varpi}}_{\Phi}^{2} \\
&\!-\!(2 \beta \delta \!-\!1) \norm{\varpi \!-\! \underline{\varpi} \!-\!(\!T_{1}\varpi \!-\! T_{1}\underline{\varpi}\!)}_{\Phi}^{2}\!.
\end{align*}
\end{lemma}
\emph{ Proof: } The proof is based on properties of  $\mathcal{A}$ and $\mathcal{B}$ in Lemma  \ref{lemma:restrictCo}, $\Phi \succ0$ by Lemma \ref{lemma:PhiEig} and resolvent properties for maximally monotone operators,  \cite{monoBookv1}, similar to Lemma 7 in \cite{PengLacraAuto2019} and Lemma 6 in \cite{Lacra} and   is omitted, due to space constraints. 
\hfill $\Box$

After establishing these properties on operators $\mathcal{A}$ and $\mathcal{B}$, 
we show next convergence to the variational GNE.

\begin{thm}\label{thmConvAggr}
Suppose Assumptions \ref{asmp:func} 
to \ref{asmp:lip} hold and ${c}\lambda_{2}\!(\!L\!) \!> \!\frac{(l^{u}_{F}\!)^{\!2}}{4\mu}$. 
Take any $\delta \! >\! \frac{1}{2\beta}$, where $\beta$ is as in Lemma  \ref{lemma:restrictCo},  $\kappa < \frac{1}{\delta}$ and  step sizes $\tau_{i}$,  $\upsilon_{i}$ and $\alpha_{i}$ chosen to satisfy Lemma \ref{lemma:PhiEig}.
	Then the action $x_{i,k}$ of each player $i\in \mathcal{N}$ generated by Algorithm 1 converges to the corresponding component in the GNE $x^*$, all agents' aggregate estimates $u_{i,k}$ converge to the same value $\sigma(x^*)$ (true aggregate) for all $i \in\mathcal{N}$, all agent's local multipliers $\lambda_{i,k}$ converge to the same multiplier $\lambda^{*}$ for all $i \in\mathcal{N}$,  corresponding  to the KKT condition  (\ref{kkt}). \end{thm}
\emph{ Proof: } Lemma 
\ref{lemma:Alg3_2_5} shows that any sequence $\{x_k, u_k, z_k, \lambda_k\}$ generated by Algorithm 1 or (\ref{eqn:AggIter})  can be obtained by some sequence  $\{x_k, u^\orth_k, z_k, \lambda_k\}$ generated by (\ref{eqn:AggIter2}), via $u_k \!=\! P_{\orthc} x_k + u^\orth_k$, cf.  (\ref{eqn:Alg5_to_3}). Thus, if we show that any sequence  $\{x_k, u^\orth_k, z_k, \lambda_k\}$ of  (\ref{eqn:AggIter2}) converges to $(\bar{x}, \bar{u}^\orth, \bar{z}, \bar{\lambda})$, then convergence  of any sequence $\{x_k, u_k, z_k, \lambda_k\}$ of Algorithm 1 follows, and moreover,  by  (\ref{eqn:Alg5_to_3}), its limit is   $(\bar{x},\! P_{\orthc} \bar{x} \!+\! \bar{u}^\orth, \bar{z}, \bar{\lambda})$, where $ P_{\orthc} \bar{x} \!=\! \mathbf{1}_N \!\otimes \! \sigma(\bar{x})$.  

To show convergence of $\{x_k, u^\orth_k, z_k, \lambda_k\}$ recall that, by Lemma  \ref{lemma:equivRep}(i), we can write any  $\varpi_k =(x_k, u^\orth_k, z_k, \lambda_k)$ of (\ref{eqn:AggIter2}) as $\varpi_{k+1} = T_{2}\circ T_{1} \varpi_{k}$, where $T_{1} \!=\! Id-\Phi^{-1}\mathcal{A}$ and $T_{2} \!=\! (Id + \Phi^{-1}\mathcal{B})^{-1}$ and, by Lemma  \ref{lemma:equivRep}(ii), any {fixed point} $(\bar{x}, \bar{u}^\orth, \bar{z}, \bar{\lambda})$ is a fixed point of $T_2 \!\circ\!T_1$ and satisfies $\bar{u}^{\orth}=\mathbf{0}_{Nn}$. By Lemma  \ref{lemma:PhiNormMono}, $T_{1}$ is restricted nonexpansive  and the inequality in Lemma  \ref{lemma:PhiNormMono}(ii) holds with respect to the consensus subspace where  $\bar{u}^{\orth}=\mathbf{0}_{Nn}$, and $T_{2} \in \mathfrak{A}(\frac{1}{2})$. Thus operators $T_1$ and $T_2$ satisfy restricted average properties as in Lemma 6 in \cite{Lacra}. Then, using an argument as in the 
proof of Theorem 25.8 in \cite{monoBookv1} for averaged operators (see 
proof of Theorem 2 in \cite{Lacra}), it follows that  any sequence $\{\varpi_k \}=\{(x_k,u^{\orth}_k,z_k,\lambda_k)\}$ converges to $\bar{\varpi} = (\bar{x},\bar{u}^{\orth},\bar{z},\bar{\lambda})$ a fixed-point of $T_2 \circ T_1$, which by Lemma \ref{lemma:equivRep}(ii),  satisfies $\bar{x}\!=\!x^*$,  $\bar{u}^{\orth}\!=\!\mathbf{0}_{Nn}$, $\bar{\lambda}\!=\!\mathbf{1}_N \!\otimes \! \lambda^*$, where $x^*$ is the variational GNE and $\lambda^*$ the corresponding multiplier.
By Lemma \ref{lemma:Alg3_2_5}, any $\{x_k, u_k, z_k, \lambda_k\}$ generated by Algorithm 1 or (\ref{eqn:AggIter})  converges, and,    by (\ref{eqn:Alg5_to_3}), its limit is  
$(x^*,\mathbf{1}_N \!\otimes \sigma(x^*),\bar{z},\mathbf{1}_N \!\otimes \lambda^*)$, hence $\{x_k\}$ converges to $x^*$, the variational GNE. 
 \hfill $\Box$

\begin{remark}
 {{  Given a globally known ${c}$ such that ${c}\lambda_{2}\!(\!L\!) \!> \!\frac{(l^{u}_{F}\!)^{\!2}}{4\mu}$, $\delta$ (depending on $\beta$ as in Lemma \ref{lemma:restrictCo}) and $\kappa$, players can independently select step sizes $\tau_i$, $\upsilon_{i}$, $\alpha_i$ to satisfy the upper bounds in  Lemma \ref{lemma:PhiEig}. 
  }} The result in Theorem  \ref{thmConvAggr} provides conditions for distributed convergence of the GNE seeking algorithm on a single time-scale. 
This is more challenging to establish than in a semi-decentralized setting where the true aggregate is provided by a coordinator. In a distributed setting, when each agent updates action and estimates the aggregate simultaneously, effectively each has to track the changing aggregate (which depends on the other agents' changing actions) while also updating its own action. 
Herein, we achieve distributed convergence on a single time-scale by using operator-splitting techniques and balancing the lack of monotonicity of the pseudo-gradient  off the consensus subspace with monotonicity of the Laplacian. 
\end{remark}
\vspace{-0.25cm}
\section{Numerical Simulations}\label{sec:example}
In this section we consider a Nash-Cournot game over a network, as in \cite{Gadjov}\cite{NedichDistAgg}\cite{parise}, for a single market with production constraints and globally coupling market capacity constraints, 
where $N=20$,  $\Omega_{i} = [0, 10]$ $\forall i\in \mathcal{N}$, $A = 1^{T}$, $b = 20$.  The cost function for each agent is $J_i(x_i,x_{-i}) \!=\! c_i(x_i) \!-\! x_if_{i}(x)$, where $c_i(x_i) \!=\! [ 1 \!+\! 2(i\!-\!1) ] x_i$ is the production cost and $f_{i}(x) \!=\! 60 - \sigma(x) - \frac{1}{2}x_{i}$ is the demand price. The variational GNE is on the boundary at $x^* = [7.809..., 5.904..., 4, 2.095..., 0.190..., 0, \dots, 0]^T$. {
For this cost function $l^{u}_{F} = 1$, $l^{x}_{F} = 1$. We compare the results with those of the algorithm in \cite{parise}, in which each agent performs $\nu$ communication rounds to update the local multiplier, then updates  its action, followed by another $\nu$ communication rounds to update its aggregate estimate. The algorithm in  \cite{parise} converges to an $\epsilon$-GNE, which approaches the variational GNE if $\nu$ goes to infinity. For simulation we chose  { $\nu = 200$} and step size $\tau = 0.01$, which means that  there are  { $400$} communication rounds before an action update, unlike Algorithm 1 where only $2$  are needed. We compare them for a star and for a ring communication network topology, respectively. { All initial conditions are set to 0.} 
} {  For Algorithm 1 we set the weighted adjacency matrix, $w_{ij} \!=\! 1$ if $(i,j)\!\in \!\mathcal{E}$, while for the algorithm \cite{parise} the mixing matrix is $W_{mix} \!= \!I \!-\! \frac{1}{20}L$ (star graph) and $W_{mix} \!=\! I \!-\! \frac{1}{3}L$ (ring graph).  }

{  
\vspace{-0.25cm}
\subsection{Star Communication Graph}
For a star communication graph topology $\lambda_{2}(L)\! = \! 1$ and $\lambda_{N}(L) \! =\!  20$, and we can use $c\! =\! 0.5$ in Lemma 4, which yields $\mu_{\tilde{\mathcal{A}}}\! =\! 0.2$, $\beta^{-1} \!=\!528.8$,  $\delta \!>\! \frac{1}{2\beta}\! =\!264$. We set 
 $\delta\! =\! 300$, $\kappa^{-1} \! =\!  500$, for which the bounds on the step sizes in Lemma 6 are $\tau_i^{-1} \! \geq \! 1806.6$, $\upsilon_i^{-1} \! \geq \!  284.4$, $\alpha_i^{-1} \! \geq \! 285.4$. Although the step sizes can be different for each agent, for simplicity they are taken to be the same: $\tau_i^{-1} \!=\! { 2000}$, $\upsilon_i^{-1}\!=\!300\!=\!\alpha_i^{-1}$. 
The agents' production plots are shown in \ref{StarDianComm}  for Algorithm 1 and in \ref{StarPariseComm} for the algorithm in \cite{parise}, respectively, where the x-axis shows the communication iterations. In the plots above the variational GNE components are denoted by stars on the y-axis to the right. In \ref{StarErrorComm} we compare the normalized error, $\norm{x-x^{*}}/\norm{x^*} 100\%$ of agents' actions from the variational GNE for the two algorithms. We see that the algorithm in \cite{parise} converges to an $\epsilon$-GNE (if $\nu$ is taken larger, one can get arbitrary close to the true GNE). 
\begin{figure}[ht]
\centering
\begin{minipage}[ht]{1\columnwidth}
	\centerline{\includegraphics[width=5cm]{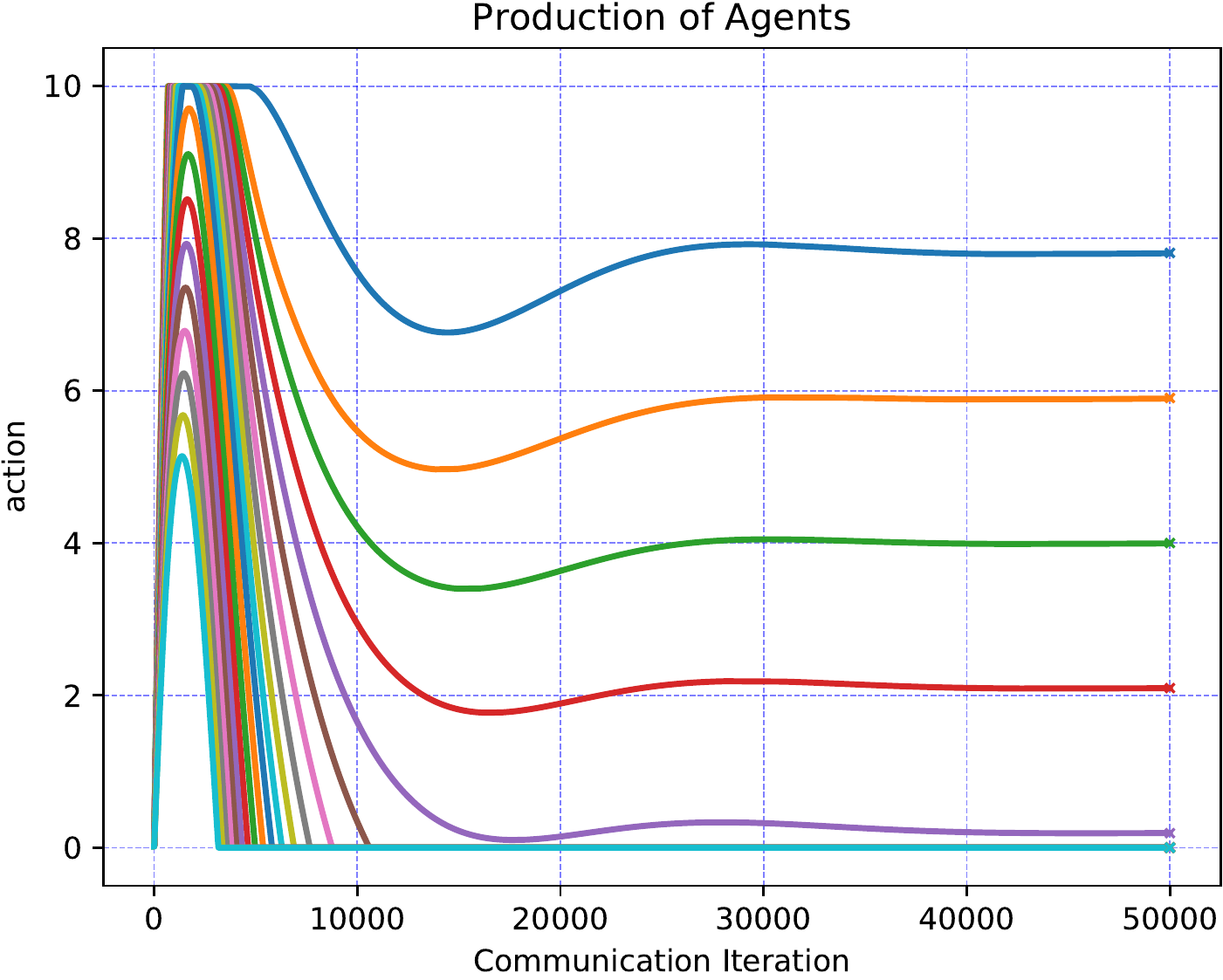}}
	\caption{Production of agents, Algorithm 1, Star topology}\label{StarDianComm}

\end{minipage}\vspace{-0.25cm}
\end{figure}
\begin{figure}[ht]
\centering
\begin{minipage}[ht]{1\columnwidth}
	\centerline{\includegraphics[width=5cm]{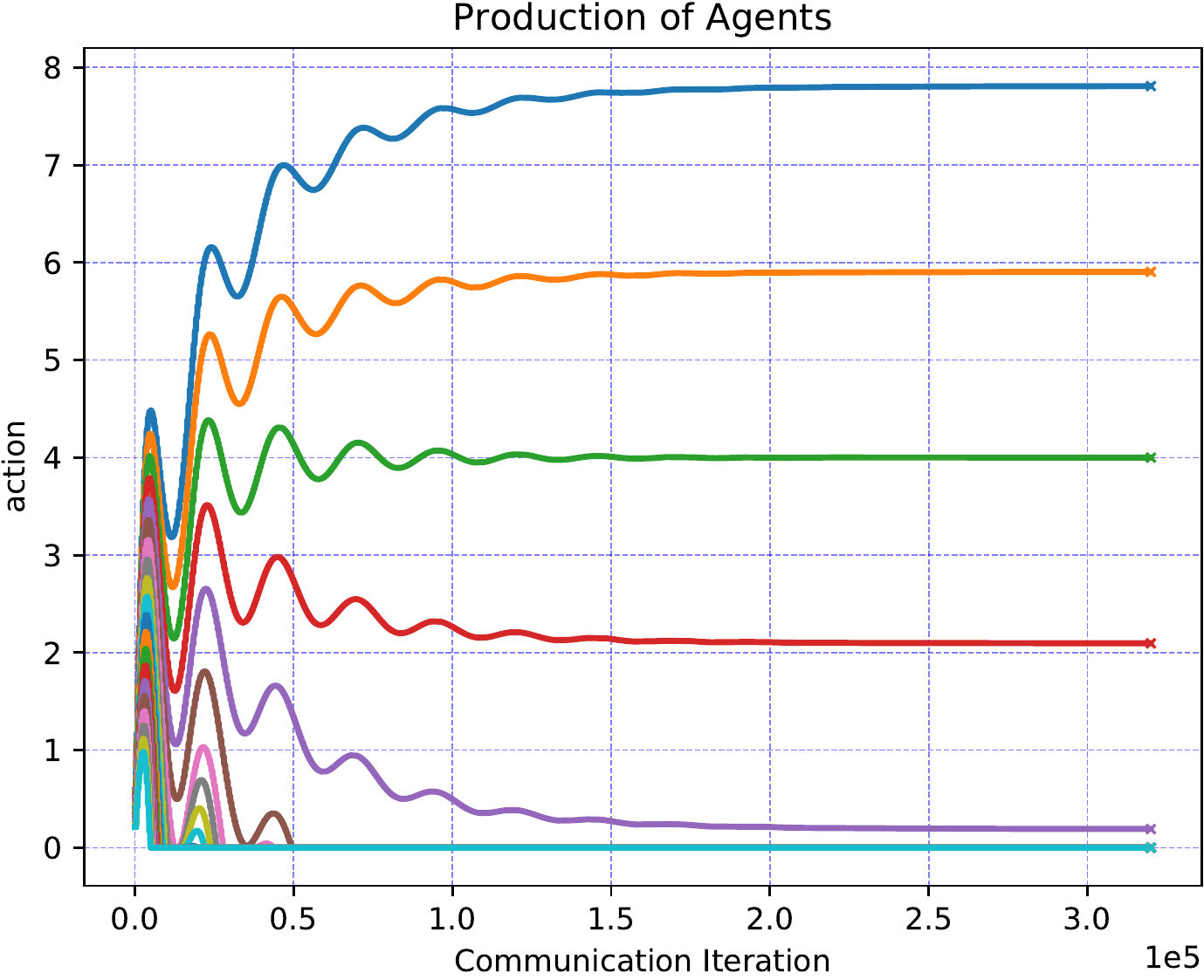}}
	\caption{Production of agents, Algorithm \cite{parise}, Star topology}\label{StarPariseComm}
\end{minipage}\vspace{-0.25cm}
\end{figure}
\begin{figure}[ht]
\centering
\begin{minipage}[ht]{1\columnwidth}
	\centerline{\includegraphics[width=5.cm]{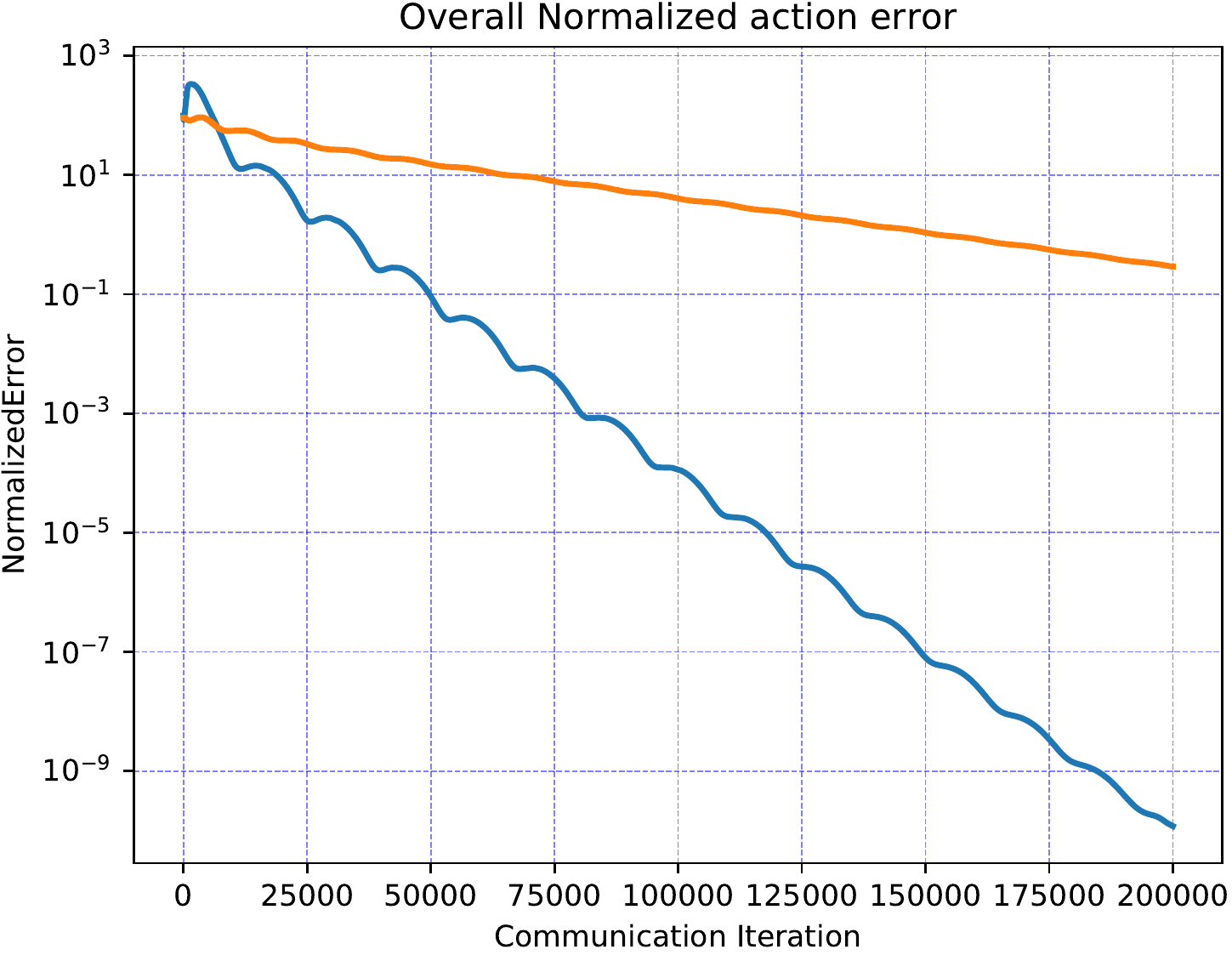}}
	\caption{Normalized error from the variational GNE, star topology. \\ \hspace*{10mm} Algorithm 1 (blue), algorithm \cite{parise} (orange)}
	\label{StarErrorComm}
\end{minipage}\vspace{-0.25cm}
\end{figure}
We note that in most cases the number of communications is the most important, but there may be situations where communications are cheap and number of action updates might be more relevant. In 
\ref{StarPariseIter} we plot  the evolution of production for  algorithm \cite{parise},  where this time the  x-axis counts the number of action updates. Thus, every action iteration  in  \ref{StarPariseIter}  corresponds to {400} communication iterations ({$\nu =200$}).  The variational GNE components are denoted by stars on the y-axis to the right. The plots for Algorithm 1 versus number of action updates look similar to those in \ref{StarDianComm} (2 communications/action update) and we omit them. 
\vspace{-0.25cm}
\begin{figure}[ht]
\centering
\begin{minipage}[ht]{1\columnwidth}
	\centerline{\includegraphics[width=5.cm]{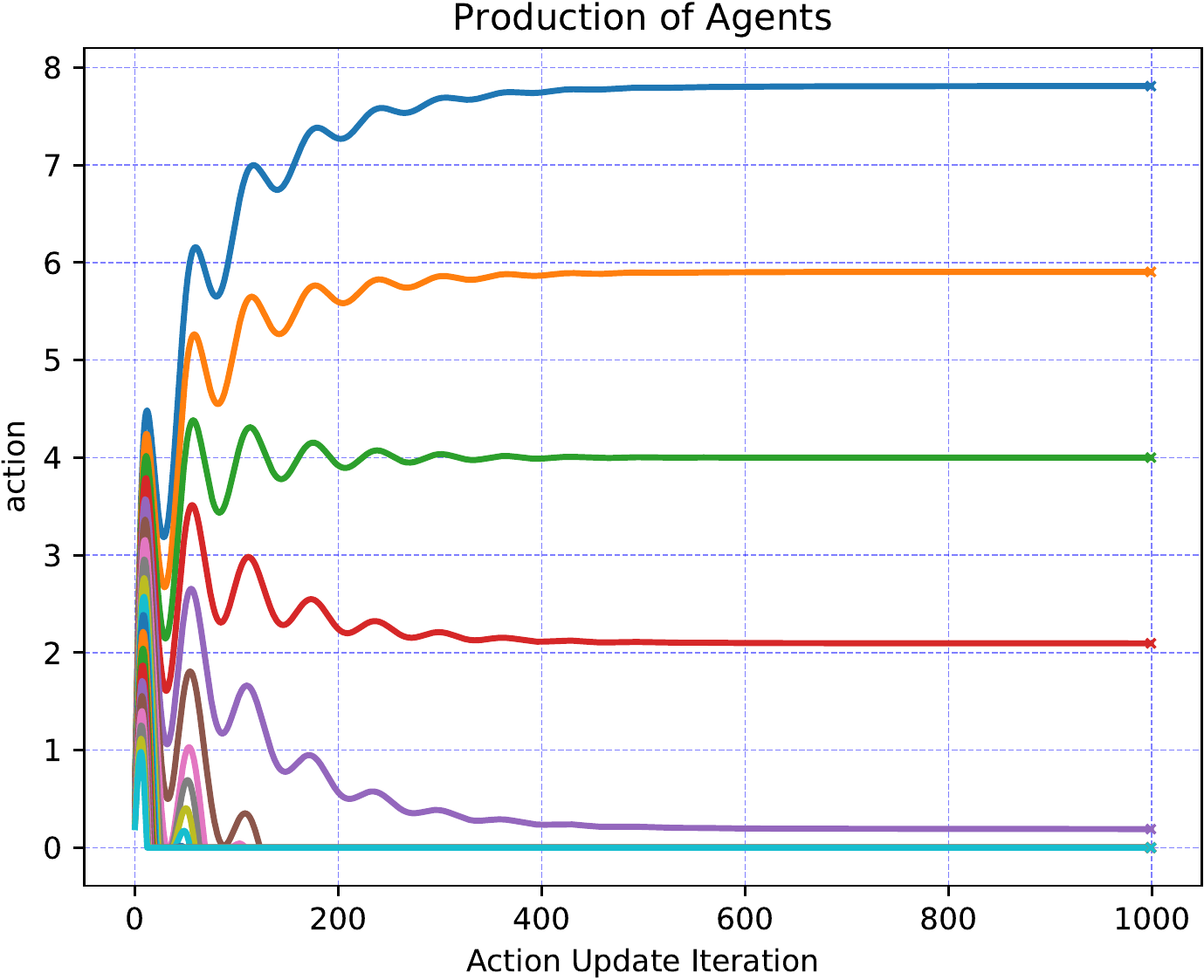}}
	\caption{Production of agents, Algorithm \cite{parise}, Star topology}\label{StarPariseIter}\vspace{-0.25cm}
\end{minipage}\vspace{-0.25cm}
\end{figure}
\vspace{-0.2cm}
\subsection{Ring Communication Graph}
For a ring communication graph  $\lambda_{2}(L)\!=\! 0.1$ and $\lambda_{N}(L) \!=\! 4$, we use $c\!=\!4$ which yields $\mu_{\tilde{\mathcal{A}}}\!=\!0.1$, $\beta^{-1} \!=\!2326$,  $\delta \!>\! \frac{1}{2\beta} \!=\!1162$. We set 
 $\delta\!=\!1200$, $\kappa^{-1} \!=\! 2000$, for which the bounds on the step sizes in Lemma 6 are $\tau_i^{-1} \!\geq  \!7942.5$, $\upsilon_i^{-1} \!\geq \! 1166.9$, $\alpha_i^{-1} \!\geq \!1167.9$. For simplicity the step sizes are taken to be the same: $\tau_i^{-1} \!=\! { 8000}$, $\upsilon_i^{-1} \!=\!1200\!=\!\alpha_i^{-1}$. 
 The agents' production plots are shown in \ref{RingDianComm} 
 for Algorithm 1 and in \ref{RingPariseComm} for the algorithm in \cite{parise}, where the x-axis counts the number of communications. 
\ref{RingErrorComm} 
compares the normalized error, $\norm{x-x^{*}}/\norm{x^*} 100\%$ from the variational GNE for the two algorithms. 
Compared to the star topology, it can be seen that convergence is slower. Thus convergence is affected by the network topology, (lower connectivity leads to slower convergence).  
On the other hand, we note that the bounds on the step sizes in Lemma 6 are conservative. In 
\ref{10_RingErrorComm} we plot the normalized error results when Algorithm 1 is run with  ten times larger step sizes, 
indicating fast convergence to the variational GNE. 
\vspace{-0.25cm}
\begin{figure}[ht]
\centering
\begin{minipage}[ht]{1\columnwidth}
	\centerline{\includegraphics[width=5.cm]{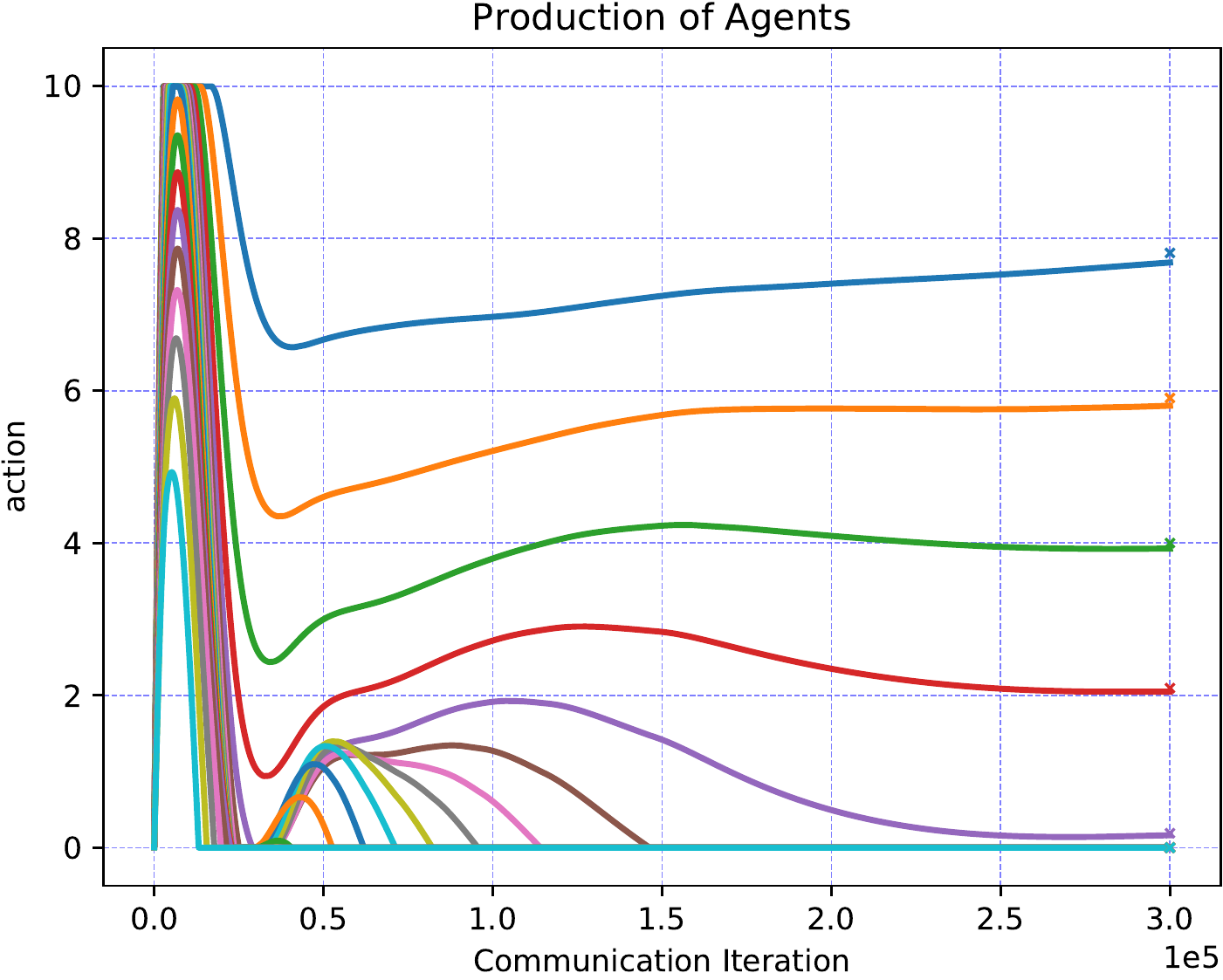}}
	\caption{Production of agents, Algorithm 1, ring topology}\label{RingDianComm}
\end{minipage}\vspace{-0.25cm}
\end{figure}
\begin{figure}[ht]
\centering
\begin{minipage}[ht]{1\columnwidth}
	\centerline{\includegraphics[width=5.cm]{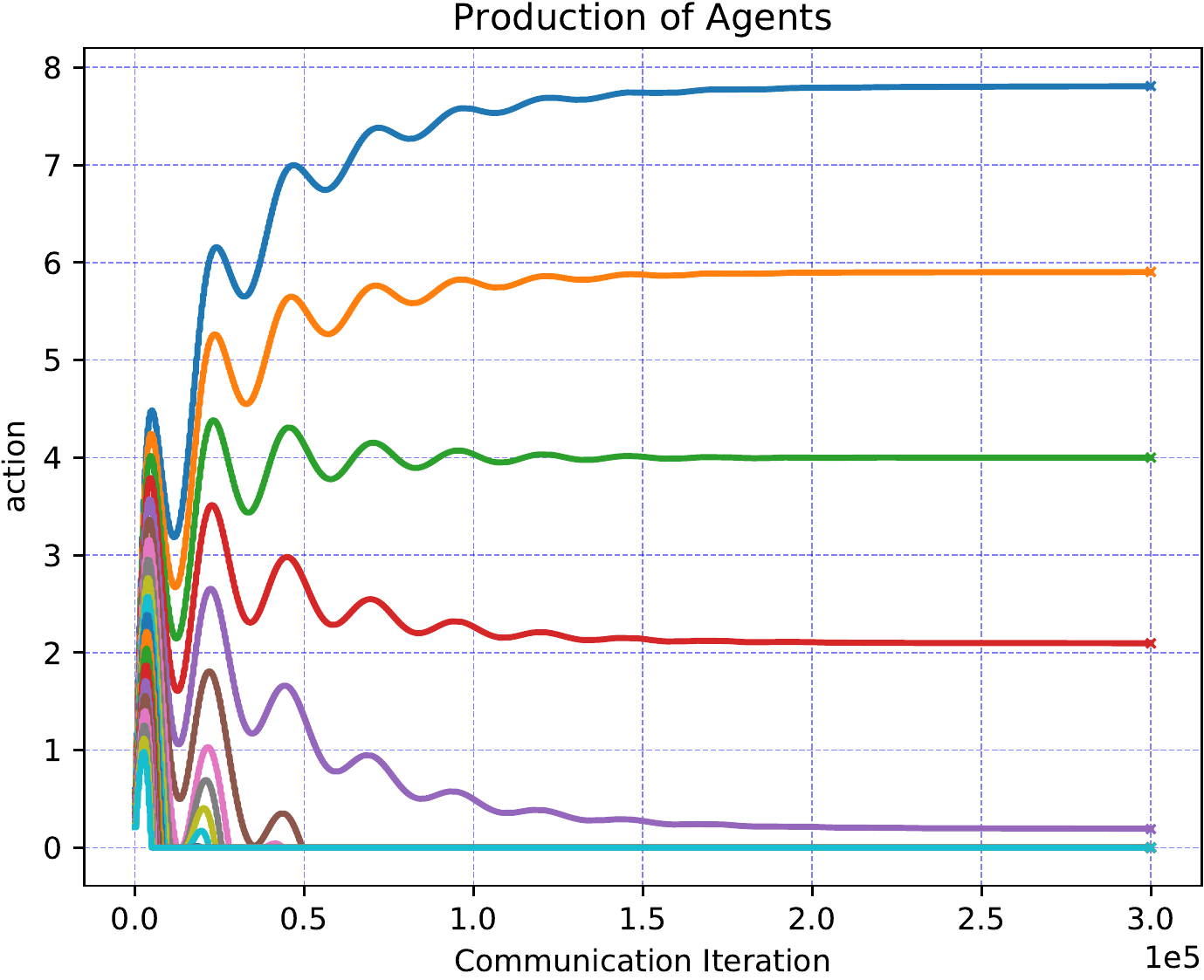}}
	\caption{Production of agents, Algorithm \cite{parise}, ring topology}\label{RingPariseComm}
\end{minipage}\vspace{-0.25cm}
\end{figure}
\begin{figure}[ht]
\centering
\begin{minipage}[ht]{1\columnwidth}
	\centerline{\includegraphics[width=5.cm]{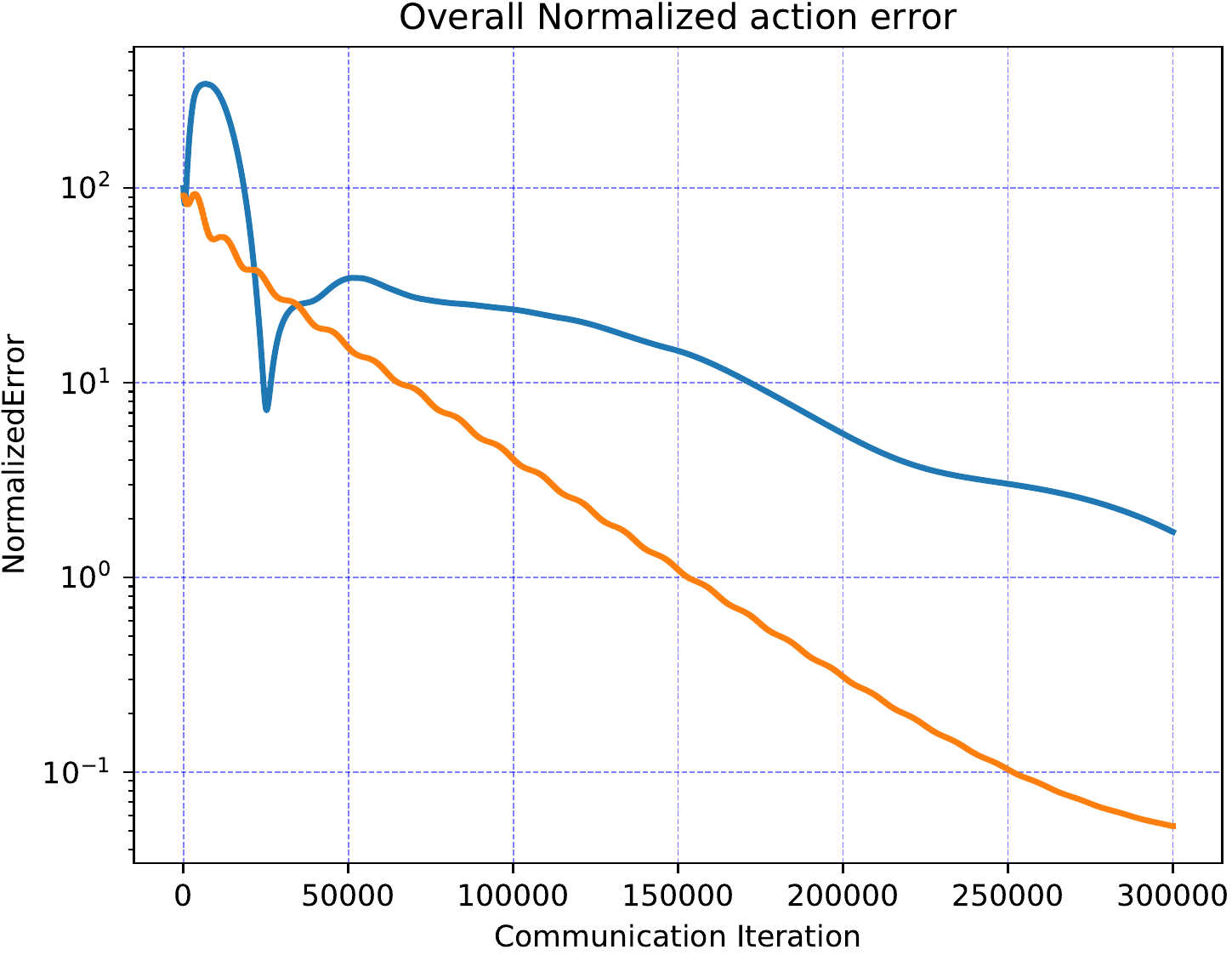}}
	\caption{Normalized error from the variational GNE, ring topology. \\ \hspace*{10mm} Algorithm 1 (blue), algorithm \cite{parise} (orange)}
	\label{RingErrorComm}
\end{minipage}\vspace{-0.25cm}
\end{figure}
\begin{figure}[ht]
\centering
\begin{minipage}[ht]{1\columnwidth}
	\centerline{\includegraphics[width=5.cm]{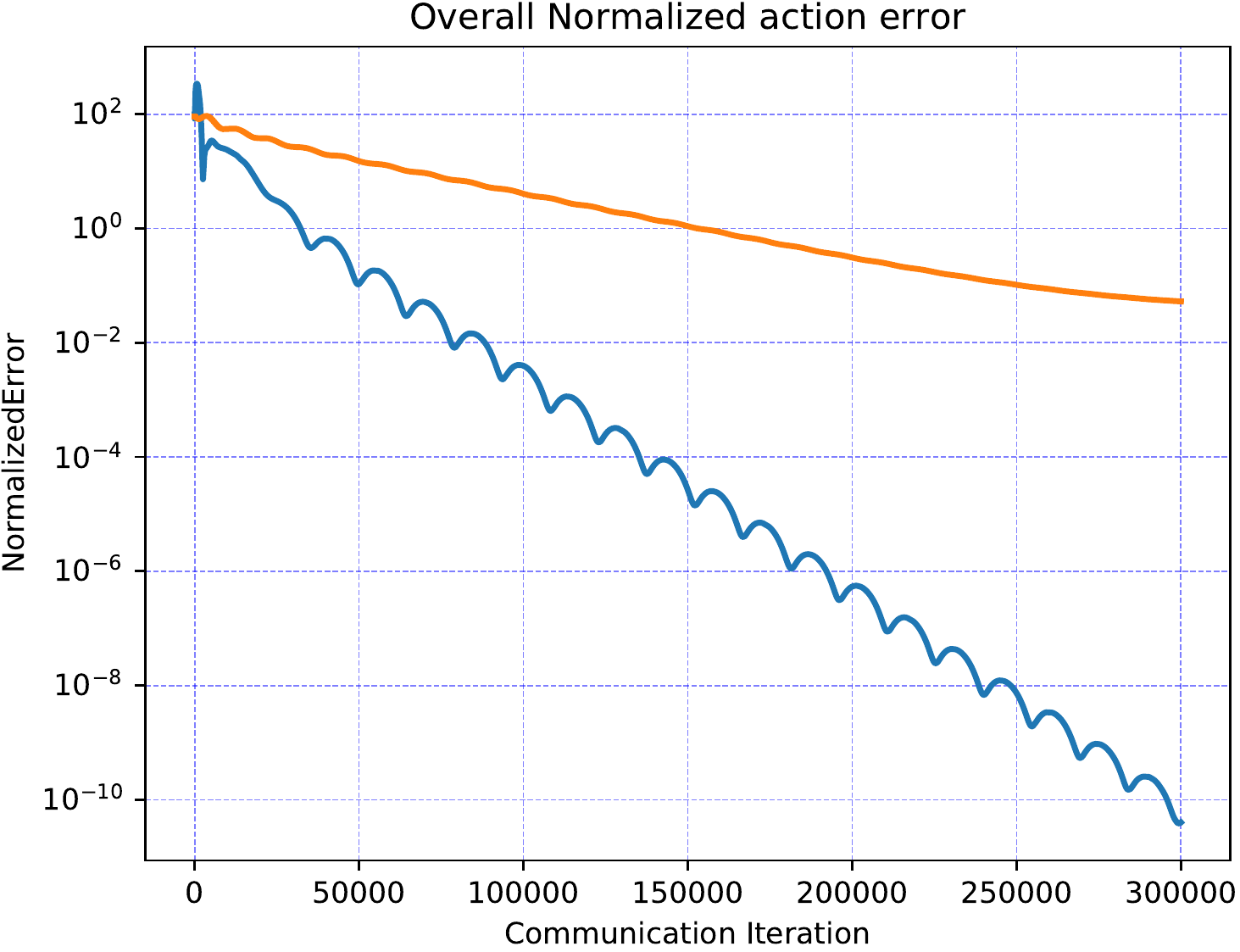}}
	\caption{Normalized error from the variational GNE, ring topology. \\ \hspace*{10mm} Algorithm 1 step sizes $\times 10$ (blue), algorithm \cite{parise} (orange)}
	\label{10_RingErrorComm}
\end{minipage}\vspace{-0.25cm}
\end{figure}

}

\vspace{-0.1cm}
\section{Conclusion}\label{sec:conclusion}
In this paper we proposed a distributed algorithm that globally converges to a variational GNE in aggregative games with affine coupling constraints. The algorithm employs  simultaneous action  and aggregate estimate update, based on local communication, and uses fixed-step sizes. We proved its convergence by forward-backward operator splitting for two preconditioned operators. We specifically exploited the invariance of the estimate average to show that the operators are restricted monotone. As future directions we can mention extension to time-varying and directed communication graphs. 
\vspace{-0.25cm}
\section*{Appendix}\label{sec:appendix}
\emph{ Proof of Lemma  \ref{lemma:Alg3_Prop}: } \hspace{-0.25cm}
(i) 	Note that from (\ref{eqn:AggIter}), using $\inp*{\mathbf{1}_N\!\otimes \!I_{n}\!}{\!\mathbf{L}_u\!} \!=\! \mathbf{0}^{T}_{Nn}$, it follows that  
$		\inp*{\mathbf{1}_N\!\otimes \!I_{n}\!}{\!u_{k\!+\!1\!}} \! = \inp*{\mathbf{1}_N\!\otimes \!I_{n}\!}{\!u_{k}\!}\!-\!\inp*{\mathbf{1}_N\!\otimes \!I_{n}\!}{\!x_{k}\!} \!+\! \inp*{\mathbf{1}_N\!\otimes  \!I_{n}\!}{\!x_{k+1}\!}, 
$ for all $k \geq 0$. Since $u_{0} \! = \!x_{0}$, the first claim follows by induction, and the second one follows by using $\sigma(u_k)\!=\!   \frac{1}{N}(\mathbf{1}^{T}_N \!\otimes \!I_{n}) u_k$, $\sigma(x_k)\!=\!   \frac{1}{N}(\mathbf{1}^{T}_N \!\otimes \!I_{n}) x_k$.

(ii) 	Let $(\overline{x}, \overline{u}, \overline{z}, \overline{\lambda})$ be a {fixed point} of Algorithm 1 or (\ref{eqn:AggIter}). Then, $	 \overline{z} \!= \! \overline{z}\! +\! \upsilon \mathbf{L}_{\lambda}  \overline{\lambda},$ 
	i.e., $\mathbf{L}_{\lambda}  \overline{\lambda} \!= \!\mathbf{0}_{Nm}$ which, by  Assumption  \ref{asmp:graph} implies that $\overline{\lambda}_{i}\!=\!\lambda^*$  for all $i \in \mathcal{N}$, for some $\lambda^*\!\in\!\reals^{m}$. From the update for $u_{k}$ in (\ref{eqn:AggIter}) if follows that 
	$ \overline{u} \!=\!  \overline{u} \!-\! \kappa {c} \mathbf{L}_{u}  \overline{u}\!+\! ({ \overline{x}\!-\! \overline{x}}), $	i.e., $\mathbf{L}_u  \overline{u} \!=\! \mathbf{0}_{Nn}$ which by Assumption  \ref{asmp:graph} implies that $ \overline{u} _i \!=\!  \underline{u}$ for all $i \in \mathcal{N}$, for some $ \underline{u} \in\reals^{n}$, hence $\sigma(\overline{u})\!=\! \underline{u}$. Using part (i) it follows  that $ \underline{u} \!=\! \sigma(\overline{x})$, i.e., in steady-state all agents have the same estimate equal to the action aggregate value. From the update of $x_{k}$ in (\ref{eqn:AggIter}) it follows that $
	\overline{x} \!=\! P_{\Omega}[ \overline{x} \!-\! \tau(\mathbf{F}(\overline{x},\overline{u}) \!+\! \Lambda^{T} \overline{\lambda})]$,  using the fact that $\mathbf{L}_u  \overline{u} \!=\! \mathbf{0}_{Nn}$. With $P_{\Omega} \!=\! (\Id \!+\! N_{\Omega})^{-1}$, this yields $ \mathbf{0}_{Nn} \! \in \! \tau^{-1} \! N_{\Omega}(\overline{x}) \!+\! \mathbf{F}(\overline{x},\overline{u}) \!+\! \Lambda^{T}\overline{\lambda}$, where   
	 $N_\Omega(x)\!=\!\Pi_{i=1}^N N_{\Omega_i}(x_i)$. Since $\tau_i\!>\!0$,  $\tau^{-1} N_{\Omega}(x)\!=\!N_\Omega(x)$, and with  $\overline{u}_i \!= \!\sigma(\overline{x})$, $\overline{\lambda}_{i}=\lambda^*$, $\forall i\in\mathcal{N}$, this yields component-wise,  
$\mathbf{0}_{n} \in \nabla_{x_i}J_{i}(\overline{x}_{i}, \sigma(\overline{x})) + A_{i}^{T}\lambda^* + N_{\Omega_{i}}(\overline{x}_{i}),\ i\in\mathcal{N}, 
$
which is the first KKT condition, (\ref{kkt}). From the update for $\lambda_{k}$ in (\ref{eqn:AggIter}) it follows that,	
$
\overline{\lambda} = P_{\reals_{+}^{Nm}}\left( \overline{\lambda}- \alpha[\mathbf{L}_{\lambda}  \overline{\lambda} + \bar{b} - \Lambda(2 \overline{x}- \overline{x}) + \mathbf{L}_{\lambda}(2 \overline{z}- \overline{z})] \right)$,  which with $\overline{\lambda}_{i}=\lambda^*$, $\forall i\in\mathcal{N}$, leads to 
$\mathbf{0}_{Nm} \in N_{\reals_{+}^{Nm}}(\mathbf{1}_N \otimes \! \lambda^*) - (\Lambda \overline{x} - \bar{b} - \mathbf{L}_{\lambda} \overline{z}), $ or $	\mathbf{0}_{Nm}=v^* - (\Lambda \overline{x} - \bar{b} - \mathbf{L}_{\lambda} \overline{z})$,  for some $v^*= col(v^*_i)_{i\in\mathcal{N}}$ with $v^*_i \in  N_{\reals_{+}^{m}}(\lambda^*)$, $\forall i\in\mathcal{N}$. Premultiplying by $ (\mathbf{1}_N^T \otimes I_m)$ and using  $(\mathbf{1}_N^T\otimes I_n) \mathbf{L}_{\lambda} = \mathbf{0}_{Nm}^{T}$ (by Assumption  \ref{asmp:graph}) and $(\mathbf{1}_N^{T}\otimes I_{m})\Lambda = A$, $(\mathbf{1}_N^{T}\otimes I_{m})\bar{b} = b$, yields $	\mathbf{0}_{m} =\sum_{i=1}^N v^*_i  -(A\overline{x}-b) $, or $\mathbf{0}_{m} \in \sum_{i=1}^N N_{\reals_{+}^{m}}(\lambda^*)\! -\! (A\overline{x} -b)$, i.e., $\mathbf{0}_{m} \! \in \! N_{\reals_{+}^{m}}(\lambda^*) \!-\! (A\overline{x} -b)$
(by Corollary 16.39 in \cite{monoBookv1}), which gives the second KKT condition, (\ref{kkt}). Therefore, $\overline{x} \!=\!x^*$  is a variational GNE and $\lambda^*$ its multiplier. 
\hfill $\Box$

{
\emph{ Proof of Lemma \ref{lemma:PhiEig} } The proof is based on the Schur complement and a diagonal dominance argument. Firstly, given any $\delta \!\!>\!\! 0$, $\kappa \!\!<\!\! \frac{1}{\delta}$, {let the step sizes} $ \tau_{i}$ be selected as in the statement, and let   $ \!  \widetilde{\tau}_{i}^{-1} \!= \! \tau_{i}^{-1}\!-  \! \delta -\!  \frac{1}{\kappa(1\!-\!\kappa\delta)} $, or 
in vector form, $\widetilde{\tau}^{-1}\!=\! \tau^{-1}  \!- 
\! (\delta +\!  \frac{1}{\kappa(1\!-\!\kappa\delta)}) I_{Nn}$, where $ \widetilde{\tau} \!=\! diag( \widetilde{\tau}_{i})_{i\in\mathcal{N}} \otimes I_{n}$. 
Using $ \tau^{-1} \!=  \!  (\delta +\!  \frac{1}{\kappa(1\!-\!\kappa\delta)}) I_{Nn}\! + \!\widetilde{\tau}^{-1}$  
in $\Phi$ (\ref{eqn:equDynOrth}) we can write 
 $\Phi \!-\!\delta I \!= \! \Phi_{1} \!+\! \Phi_{2} $ where \vspace{-0.22cm} 
\begin{align*}
\Phi_{1} &= \begin{bmatrix}
 \frac{1}{\kappa(1\!-\!\kappa\delta)}  I_{Nn} + \kappa^{-1} P_{\orth}   & - P_{\orth} \kappa^{-1} & \mathbf{0} & \mathbf{0} \\  	
-\kappa^{-1}P_{\orth}  & (\kappa^{-1}\! -\! \delta ) I_{Nn} & \mathbf{0} & \mathbf{0} \\
\mathbf{0} & \mathbf{0} & \mathbf{0} & \mathbf{0} \\
\mathbf{0} & \mathbf{0} & \mathbf{0} & \mathbf{0}
\end{bmatrix} 
\end{align*}\vspace{-0.35cm}
\begin{align*}
\Phi_{2} &= \begin{bmatrix} \,
	 \widetilde{\tau}^{-1} & \mathbf{0} & \mathbf{0} & -\Lambda^T \\
	\mathbf{0} & \mathbf{0} & \mathbf{0} & \mathbf{0} \\
	\mathbf{0} & \mathbf{0} & \upsilon^{-1} - \delta I_{Nm} & \mathbf{L}_{\lambda} \\
	-\Lambda & \mathbf{0} & \mathbf{L}_{\lambda} & \alpha^{-1} - \delta	 I_{Nm}
	\end{bmatrix}
\end{align*}
and $\mathbf{0}$ is an all-zero matrix of appropriate dimension. 

Since  $\kappa^{-1} \!-\! \delta \!>\!0$, using Schur's complement, 
the top-left block matrix in $\Phi_1$ is positive semi-definite if  \vspace{-0.25cm}
$$
 \frac{1}{\kappa(1\!-\!\kappa\delta)}  I_{Nn} + \kappa^{-1} P_{\orth}  \! -\! P_{\orth} \kappa^{-1}(\kappa^{-1}\! - \!\delta)^{-1}\kappa^{-1}P_{\orth}  \!\succeq \!0
$$ 
or, 
equivalently, if 
$
 \kappa^{-1} P_{\orth}  \!+ \frac{1}{\kappa(1 \!-\!\kappa\delta)} (I_{Nn}\! -\!  P_{\orth} P_{\orth} )\succeq 0, 
$
which holds since $P_{\orth}$ is a projection matrix. 
Thus $\Phi_{1}$ is positive semi-definite.  

By Ghershgorin's theorem, 
a sufficient condition for $\Phi_2$ to be positive semi-definite is to be diagonally dominant, i.e., 
\vspace{-0.25cm}
\begin{equation}
\begin{array}{lll}\label{equ_lem5_4_1}
 \widetilde{\tau}_{i}^{-1}   & \geq  &   \max_{j=1,\dots,n} \{\sum^m_{k=1} | [A_i^T]_{jk} |\}  \\
\nu_i^{-1}- \delta   & \geq  &   \sum_{j=1}^m  | L_{ij} | = 2 d_i \\
\sigma_i^{-1}-\delta &  \geq &   \max_{j=1,\dots,m} \{ \sum^{n}_{k=1} | [A_i]_{jk} | \} + \sum_{j=1}^m  | L_{ij} | \\
\qquad                    & = &\max_{j=1,\dots,m} \{ \sum^{n}_{k=1} | [A_i]_{jk} | \}+2d_i
\end{array}
\end{equation}
for all $ i \!\in \! \mathcal{N}$. 
Since 
  $\widetilde{\tau}_{i}^{-1} \! = \! \tau_{i}^{-1}\!-
\! \delta  -\!  \frac{1}{\kappa(1\!-\!\kappa\delta)} $,  if $\tau_i$ $\nu_i$ and $\alpha_i$ are selected as in the lemma statement, it follows  that 
all inequalities in (\ref{equ_lem5_4_1})  hold, hence 
$\Phi_{2}$ is positive semi-definite. 
Thus, $ \Phi -\delta I$ is  positive semi-definite, and therefore $\Phi$ is  positive definite. \hfill $\Box$
}

\vspace{-0.25cm}
\bibliographystyle{IEEEtran}
\bibliography{referencesCDC}

\end{document}